\documentclass[11pt]{amsart}

\usepackage{amssymb}
\usepackage{amsmath}
\usepackage{amsthm}
\usepackage{amsfonts}

\newcommand{\A}{{\mathcal{A}}}
\newcommand{\B}{{\mathcal{B}}}
\newcommand{\D}{{\mathcal{D}}}
\newcommand{\E}{{\mathbb{E}}}

\newcommand{\F}{{\mathcal{F}}}

\newcommand{\K}{{\mathcal{K}}}
\newcommand{\kk}{\kappa}
\newcommand{\Ha}{{\mathcal{H}}}

\newcommand{\C}{{\mathbb{C}}}
\newcommand{\R}{{\mathbb{R}}}
\newcommand{\T}{{\mathbb{T}}}
\newcommand{\M}{{\mathcal{M}}}
\newcommand{\NN}{{\mathcal{N}}}
\newcommand{\N}{{\mathbb{N}}}

\newcommand{\gwia}{^{\star}}
\newcommand{\jed}{{\mathbf{1}}}

\DeclareMathOperator{\NC}{NC}
\DeclareMathOperator{\spann}{span}
\DeclareMathOperator{\tr}{tr}

\DeclareMathOperator{\Alg}{Alg}

\theoremstyle{plain}
\newtheorem{lemma}{Lemma}
\newtheorem{theorem}[lemma]{Theorem}
\newtheorem{thedef}[lemma]{Theorem and Definition}

\newtheorem{corollary}[lemma]{Corollary}
\newtheorem{conjecture}[lemma]{Conjecture}

\theoremstyle{definition}
\newtheorem*{definition}{Definition}

\theoremstyle{remark}
\newtheorem*{remark}{Remark}

\begin{document}

\author{ Piotr \'Sniady}
\address{Institute of Mathematics,
University of Wroclaw, pl.\ Grunwaldzki 2/4, 50-384 Wroclaw, Poland}
\email{Piotr.Sniady@math.uni.wroc.pl}

\author{Roland Speicher}
\address{Department of Mathematics and Statistics, Queens University,
Kingston Ontario K7L 3N6, Canada}
\email{speicher@mast.queensu.ca}

\title[Continuous Family of Invariant Subspaces]
{Continuous Family of Invariant Subspaces for $R$--diagonal Operators}

\begin{abstract}
We show that every $R$--diagonal operator $x$ has a continuous family of
invariant subspaces relative to the von Neumann algebra generated by $x$.
This allows us to find the Brown measure of $x$ and to find a new
conceptual proof that Voiculescu's $S$--transform is multiplicative.
Our considerations base on a new concept of $R$--diagonality with
amalgamation, for which we give several equivalent characterizations.
\end{abstract}

\maketitle

\section {Introduction}
\label{chap:introduction}
\subsection {Upper triangular matrices}
\subsubsection {Finite matrices}
It is a not very difficult exercise to show that every matrix
$x\in\M_n(\C)$ is unitarily equivalent to an upper triangular matrix
\index{upper triangular matrix} $y=y(x)$
of the form
\begin{equation} y= \left[ \begin{array}{ccccc}
\lambda_1 & a_{1,2} &   \cdots & a_{1,n-1} & a_{1,n} \\
0 & \lambda_2 & \cdots & a_{2,n-1} & a_{2,n} \\
\vdots&                 &       \ddots & \vdots   & \vdots        \\
& & & \lambda_{n-1} & a_{n-1,n} \\
0&       & \cdots &             0               & \lambda_n
\end{array} \right] ,
\label{eq:utm}
\end{equation}
where $\lambda_1,\dots, \lambda_n$ are the eigenvalues of $x$
and $a_{i,j}\in\C$ ($i<j$). By $(e_1,\dots,e_n)$ we denote the orthonormal
base of $\C^n$ in which the matrix $x$ takes the form
(\ref{eq:utm}).

This upper triangular representation is usually not unique due to the
possibility of changing the order of eigenvalues. Now we shall choose the
order of eigenvalues so that $|\lambda_1|<|\lambda_2|<\cdots<|\lambda_n|$
(for simplicity we shall not consider the case when some absolute values of
eigenvalues are equal); let corresponding eigenvectors be
$v_1,\dots,v_n$.
For $s\geq 0$ let us consider the set
\begin{equation} V_s=\left\{v\in \C^n:
\lim_{k\rightarrow\infty} \sqrt[k]{\| x^k v\|}\leq s \right\} .
\label{eq:niezmienniczaC}
\end{equation}
One can show
that $V_s$ is a linear space and that the
orthonormal base $(e_1,\dots,e_n)$, the base of eigenvectors
$(v_1,\dots,v_n)$, and the spaces $V_s$ are closely related by the following
property:
if $|\lambda_i|\leq s<|\lambda_{i+1}|$ then $V_s=\spann\{e_1,\dots,
e_i\}=\spann\{v_1,\dots,v_i\}$.

\subsubsection {Upper triangular form property}
Much information can be obtained directly from the upper triangular form
(\ref{eq:utm}) of a matrix $x$, for example its spectrum and relations
between the matrix and its adjoint. Therefore it is natural to ask if it is
possible to find some kind of the upper triangular form for operators acting
on general Hilbert spaces.
For every $s\geq 0$ and $x\in B(\Ha)$ we can define an analogue of
the linear space (\ref{eq:niezmienniczaC}): \begin{equation} V_s=\overline{
\left\{v\in \Ha: \limsup_{k\rightarrow\infty} \sqrt[k]{\| x^k v\|}\leq s
\right\}} . \label{eq:niezmienniczaHa}
\end{equation}
These sets $V_s$ are always closed linear spaces. They
were introduced and investigated by Dykema and Haagerup
\cite{DyHa2000}, who proved that $V_s$ is always
an invariant subspace of the operator $x$. Let $p_s$
denote the orthogonal projection from $\Ha$ onto $V_s$. They also showed that
if $x$ is an element of a von Neumann algebra $\NN\subseteq B(\Ha)$ then 
$p_s\in\NN$ and that $p_s$ regarded as an abstract element of $\NN$ does not
depend on the representation of the algebra $\NN$.

If $0\leq s_1< s_2<\cdots<s_{n-1}$ then
$V_{s_1}\subseteq V_{s_2}\subseteq\cdots\subseteq V_{s_{n-1}}$ are invariant
subspaces of the operator $x$, therefore $x$ can be written as
\begin{equation}
\label{eq:utm1}
x=\left[ \begin{array}{ccccc} e_1 x e_1 &
e_1 x e_2 & \cdots & e_1 x e_{n-1}& e_1 x e_n \\ 0 & e_2 x e_2 &\cdots & e_2
x e_{n-1}& e_2 x e_n \\ \vdots& & \ddots & \vdots & \vdots \\ & & & e_{n-1} x
e_{n-1} & e_{n-1} x e_n \\ 0& & \cdots & 0 & e_n x e_n \end{array} \right], 
\end{equation}
where $e_1=p_{s_1}$, $e_2=p_{s_2}-p_{s_1}$,
$e_3=p_{s_3}-p_{s_2}$, \dots, $e_n=I-p_{s_{n-1}}$ are orthogonal projections 
with the property that $e_1+\cdots+e_n=I$.
In the finite dimensional case the diagonal elements
$\lambda_1,\dots,\lambda_n$ of the matrix (\ref{eq:utm}) are ordered in such
a way that their absolute values increase, therefore we would expect that
in the general case the following analogue holds: 
\begin{equation}
\sigma_i\subseteq \{z\in\C: s_{i-1} \leq |z| \leq s_{i} \} \qquad\mbox{for
every } 1\leq i\leq n,
\label{eq:spektra2}
\end{equation} 
where $\sigma_i$ denotes the spectrum of $e_i x e_i$ regarded as an operator
acting on the space $e_i \Ha$ and $s_0=0$, $s_n=\infty$.

If an operator $x$ has the property that for every $0<s_1<s_2<\dots<s_{n-1}$
inclusions (\ref{eq:spektra2}) hold, we say that $x$ has the {\bf upper
triangular form property} \index{upper triangular form property}.
 
It is not difficult to find examples of operators which do not
have this property.

\subsection {Applications}
Operators with the upper triangular form property have a very nice structure
and
therefore many questions, for example the invariant subspace problem and the
computation of the Brown measure, are for these operators particularly easy.

\subsubsection {Invariant subspace problem}
The invariant subspace problem relative to a von Neumann algebra
$\NN\subseteq B(\Ha)$ asks if every operator $x\in\NN$ has a proper,
nontrivial invariant subspace $\Ha_0\subseteq \Ha$ such that the orthogonal
projection $p$ onto $\Ha_0$ is an element of $\NN$.
Originally this problem was formulated for the case $\NN=B(\Ha)$ and in
both formulations it remains open.

It was pointed out by Dykema and Haagerup \cite{DyHa2000} that if an
operator $x$ has the upper triangular form property and if we can choose $s$
such that $\inf_{\lambda\in\sigma(x)}
|\lambda|<s<\sup_{\lambda\in\sigma(x)} |\lambda|$
(where $\sigma(x)$ denotes the spectrum of $x$) then the space $V_s$ is
the wanted nontrivial invariant subspace.

\subsubsection {Brown measure}
\index{Brown measure}
Suppose that $x$ is an element of a finite von Neumann algebra $\NN$ with a
normal faithful tracial state $\phi$. The Fuglede--Kadison
determinant \index{Fuglede--Kadison determinant} (cf \cite{FK}) defined by
$$\Delta(x)=\exp\big[\phi \big(\ln (x x\gwia)^{\frac{1}{2}} \big) \big] $$
is a generalization of a determinant of a finite matrix.

The Brown measure
$\mu_x$ of the element $x$ is a Schvarz distribution on the complex plane
defined by (cf \cite{Brown,HL})
$$\mu_x= \frac{1}{2\pi} \left(\frac{\partial^2}{\partial a^2}+
\frac{\partial^2}{\partial b^2} \right) \ln \Delta[x-(a+i b)I ].$$
One can show that in fact $\mu_x$ is a
probability measure supported on a subset of the spectrum of $x$.
It can be regarded as a generalization of the normalized counting measure
$\frac{1}{n}
(\delta_{\lambda_1}+\delta_{\lambda_2}+\cdots+\delta_{\lambda_n})$ that is
assigned to a finite matrix with eigenvalues
$\lambda_1,\lambda_2,\dots,\lambda_n$.

\subsubsection {Brown measure for upper triangular matrices}
\label{subsect:miarabrowna}
For a matrix $y\in \M_n(\C)$ in the form (\ref{eq:utm}) we have
$\det y=\lambda_1 \lambda_2 \cdots \lambda_n$; similarly
the Fuglede--Kadison determinant of $x\in\NN$ in the form (\ref{eq:utm1})
is a weighted product of determinants of diagonal submatrices:
$$\Delta(x)=[\Delta_{e_1} (e_1
x e_1)]^{\phi(e_1)} [\Delta_{e_2} (e_2 x e_2)]^{\phi(e_2)} \cdots
[\Delta_{e_n} (e_n x e_n)]^{\phi(e_n)},$$
where $\Delta_{e_k} (e_k x e_k)$
denotes the Fuglede--Kadison determinant of $e_k x e_k$ regarded as an
element of the algebra $e_k \NN e_k$ with state $\frac{1}{\phi(e_k)} \phi$
(cf Proposition 1.8 in \cite{Brown}).

Therefore we have that the Brown measure of $x$ is a weighted sum of Brown
measures of diagonal elements:
\begin{equation} \mu_x=\phi(e_1) \mu_{e_1; e_1 x e_1}+\phi(e_2)
\mu_{e_2; e_2 x e_2}+\cdots+\phi(e_n) \mu_{e_n; e_n x e_n},
\label{eq:miarabrowna1}
\end{equation}
where $\mu_{e_k; e_k x e_k}$ denotes the Brown measure of $e_k x e_k$
regarded as an element of the algebra $e_k \NN e_k$ with state
$\frac{1}{\phi(e_k)} \phi$ (cf Theorem 4.3 in \cite{Brown}).
If furthermore condition (\ref{eq:spektra2}) for the spectra of diagonal
elements holds, we have that the supports of the measures $\mu_{e_k; e_k x e_k}$
($1\leq k\leq n$) have disjoint interiors and
\begin{equation} \mbox{supp}(\mu_{e_k; e_k x e_k})\subseteq\{ z\in\C :
s_{k-1}\leq |z|\leq s_k\}.
\label{eq:miarabrowna2}
\end{equation}

The above properties (\ref{eq:miarabrowna1}) and (\ref{eq:miarabrowna2})
provide complete information on the radial part of the measure $\mu_x$;
if we know that $\mu_x$ is rotation invariant, they uniquely
determine $\mu_x$ to be the product measure $\nu_x \times \chi$ on
$[0,\infty)\times {\mathbb{T}}=\C$ (these sets are identified by the polar
coordinates $(r,\psi)\mapsto r e^{i\psi}$), where $\chi$ is the uniform
measure on ${\mathbb{T}}$ and $\nu_x$ is the
probability measure on $[0,\infty)$ with the distribution function
$$F_x(s)=\nu_x[0,s]=\phi(p_s)\qquad\mbox{for every } s\geq 0. $$

\subsection {Results and methods}
The above considerations show how important it is to characterize the class
of operators with the upper triangular form property. Dykema and Haagerup (cf
\cite{DyHa2000}) showed that the circular element and more generally every
free circular Poisson element has the upper triangular form property. Their
proof was based on explicit random matrix models for these operators. In
this article we show a significantly simpler proof of a more general
result that every $R$--diagonal operator has the upper triangular form
property.


$R$--diagonal operators were introduced by Nica and Speicher
\cite{NiSp1997A} in the setup of the free probability theory of Voiculescu
\cite{VDN}; we say that $x\in\NN$ is $R$--diagonal with respect to
$\phi$ (where $\NN$ is a finite von Neumann algebra with a trace $\phi$) if
the distribution of $x$ is invariant under multiplication with free
unitaries. Therefore a typical $R$--diagonal element generates the
$\mbox{II}_1$ factor $L({\mathbb{F}}_2)$ associated to the nonabelian free
group ${\mathbb{F}}_2$ with two generators.

Our considerations use the operator--valued version of
free probability. The key idea is to relate properties of
noncommutative random variables regarded as elements of a scalar probability
space with corresponding properties of these random variables regarded as
elements of an operator--valued probability space.

In Section  \ref{chap:preliminaries} we recall shortly some basic concepts of
the scalar and operator--valued free probability theory.

In Section  \ref{chap:rdiagonality}
we define a new concept of the $R$--diagonality with
amalgamation, which is a generalization of the $R$--diagonality.
In Theorem and Definition \ref{thedef:centralne}, which is a generalization
of the result of Nica, Shlyakhtenko and Speicher \cite{NiShSp}, we show
a few equivalent conditions on the $R$--diagonality in the operator--valued
case. In Theorem \ref{theo:nowardiagonalnosc} we present a new
characterization of the scalar $R$--diagonality in terms of freeness.

Section  \ref{chap:compatible} is devoted to the problem how to
relate properties of noncommutative random variables regarded as elements
of a scalar probability space with corresponding properties of these random
variables regarded as elements of an operator--valued probability space.
We recall some old results in this field and show
Theorem \ref{theo:rdiagonalnosc}, one of the main results of this article,
that if $x$ is $R$--diagonal with amalgamation over certain
algebra $\D$ and if $\{x,x\gwia\}$ and the algebra $\D$ are free then $x$ is
scalar $R$--diagonal as well. As an immediate application
we discuss the upper triangular representation of circular free Poisson
elements which was found by Dykema and Haagerup \cite{DyHa2000}.
We show also a useful Lemma \ref{lem:pierwiastek} that every positive
operator has a square root which takes the upper triangular form with
respect to a given family of orthogonal projections.

In Section  \ref{chap:main} we
prove the central result of this article that every $R$--diagonal operator
has the upper triangular form property and we find a new way to
obtain a result of Haagerup and Larsen \cite{HL,La}: the explicit formula
for the Brown measure of $R$--diagonal elements.

In Section
\ref{chap:produkt} we relate the upper triangular forms for
$x$, $y$ and $xy$, where $x$ and $y$ are free $R$--diagonal operators and
we present a new interesting proof of the well known fact that Voiculescu's
$S$--transform is multiplicative.

Finally, in Section  \ref{chap:losowemacierze}
we show a possible direction of the future research: the problem how to
relate the properties of $R$--diagonal elements and their random matrix models.


\section {Preliminaries}
\label{chap:preliminaries}
\subsection [Probability spaces]{Scalar and operator--valued probability spaces}
\subsubsection {Basic definitions}
In the following we shall consider only unital algebras over $\C$.
We shall also consider only unital inclusions of unital algebras
$\A_1\subseteq\A_2$, i.e.\ the identity $I$ of the algebra $\A_1$ is also the
identity of the algebra $\A_2$. Identities of all algebras will be denoted
by $I$.

A scalar probability space \index{scalar probability space}
(also called a noncommutative probability space)
\index{noncommutative probability space}
is a pair $(\A,\phi)$ consisting of a unital $\star$--algebra $\A$ and a
normal, tracial, faithful state $\phi:\A\rightarrow\C$.

A pair $(\D\subseteq\A,E)$ is an operator--valued probability space
\index{operator--valued probability space} (or
a $\D$--valued probability space) if $\A$ is a unital $\star$-algebra,
$\D\subseteq\A$ is a unital $\star$--subalgebra, and $E:\A\rightarrow\D$ is a
conditional expectation (i.e.\ $E$ is linear, satisfies $E(I)=I$ and
$E(d_1 x d_2)=d_1 E(x) d_2$ for every $d_1,d_2\in\D$ and $x\in\A$).

The concept of the operator--valued probability space is a
generalization of the scalar probability space in which the elements of
the algebra $\D$ play the role of scalars. In particular, every scalar
probability space $(\A,\phi)$ can be regarded as
a $\C$--valued probability space $(\C\subseteq\A,\phi)$, but
not every $\C$--valued probability space is a scalar probability
space.

An element $a\in\A$ will be referred to as a noncommutative random
variable; \index{noncommutative random variable}
$\phi(a)$ and $E(a)$ will be called the expectation of $a$.
\index{expectation, expected value}

\subsubsection {Notation}
For simplicity we will use the following notation: by $x\in (\D\subseteq\A,E)$
we shall mean
that $x\in\A$ and that the noncommutative random variable $x$ is
regarded as an
element of a probability space $(\D\subseteq\A,E)$;
similarly for $x\in (\A,\phi)$.

Also, by the inverse (or spectrum, etc.) of $x\in\A$ we shall mean the inverse
(resp.\ spectrum, etc.) of $x$ regarded as an element of a unital algebra $\A$.

The spectrum of the operator $x$ will be denoted by $\sigma(x)$.

\subsubsection {Enlargement of a $\D$--valued probability space}
The $\D$--valued probability space $(\D\subseteq\tilde{\A},\tilde{E})$ is an
enlargement of the $\D$--valued probability space \index{enlargement of an
operator--valued probability space} $(\D\subseteq\A,E)$ if
$\A\subseteq\tilde{\A}$ and the restriction of the expectation value
$\tilde{E}:\tilde{\A}\rightarrow\D$ to the algebra $\A$ coincides with the
expectation $E:\A\rightarrow\D$.

\subsubsection {Operator--valued distribution}
\index{operator--valued distribution}
We say that elements $x_1\in(\D\subseteq\A_1,E_1)$ and
$x_2\in(\D\subseteq\A_2,E_2)$ are
identically distributed if for every $n\geq 1$,
all $s_1,\dots,s_n\in\{1,\star\}$, and all $d_1,\dots,d_{n-1}\in\D$ we have
$$E_1(x_1^{s_1} d_1 x_1^{s_2} d_2 \cdots d_{n-1} x_1^{s_n})=
E_2(x_2^{s_1} d_1 x_2^{s_2} d_2 \cdots d_{n-1} x_2^{s_n}).$$

\subsection {Basics of the free probability theory}
\subsubsection {Freeness and freeness with amalgamation}
In this section we recall briefly basics of the operator--valued
free probability theory
\cite{Voi1985,VoiAst,Sp1998,VDN}.

Let $(\A_s)_{s\in S}$ be a family of subalgebras of $\A$ such that
$\D\subseteq\A_s$ for every $s\in S$. We say that the algebras
$(\A_s)_{s\in S}$ are free with amalgamation over $\D$ (or simply
free over $\D$)
\index{freeness with amalgamation} if:
\begin{equation} E(x_1 x_2 \cdots x_n)=0  \label{eq:freeness}
\end{equation}
holds for every $n\geq 1$, every $s_1,s_2,\dots,s_n\in S$ such
that $s_1\neq s_2$, $s_2\neq s_3$,\dots, $s_{n-1}\neq s_n$, and every
$x_1\in\A_{s_1},\dots, x_n\in\A_{s_n}$ such that
$E(x_1)=\cdots=E(x_n)=0$ (cf \cite{VDN,VoiAst}).

Let $(X_s)_{s\in S}$ be a family of subsets of $\A$. We say that the sets
$(X_s)_{s\in S}$ are free with amalgamation over $\D$ if the algebras
$\A_s:=\mbox{Alg}(X_s\cup \D), s\in S$, are so.

Let $(\A,\phi)$ be a scalar probability space.
If a family $(X_s)_{s\in S}$ of subsets of
$\A$ is free with amalgamation over $\C$, we simply say that it is free.
\index{freeness}

\subsubsection {Noncrossing partitions. Complementary partition}
If $X$ is a finite, ordered set, we denote by $\NC(X)$ the set of all
noncrossing partitions of $X$ \index{noncrossing partitions}
(cf \cite{Sp1998, Kreweras}). A noncrossing partition $\pi=\{Y_1,\dots,Y_n\}$
of $X$ is a decomposition of $X$ into disjoint sets:
$$X=Y_1\cup\cdots\cup Y_n, \qquad Y_i\cap Y_j=
\emptyset \quad \mbox{ if }i\neq j,$$
which has the additional property that
if $a<b<c<d\in X$ are such that $a,c\in Y_i$ and $b,d\in Y_j$ then $i=j$.

In the set $\NC(X)$ we introduce a partial order: for
$\pi,\sigma\in\NC(X)$ we say that $\pi\prec\sigma$ iff for every $A\in\pi$
there exists $B\in\sigma$ such that $A\subseteq B$.

If $\sigma,\tau\in\NC(X)$ then by $\sigma\vee\tau$ we denote the (unique)
minimal element in the set $\{\pi\in\NC(X): \sigma\prec\pi \mbox{ and }
\tau\prec\pi\}$.

By ${\bf 1}_X\in\NC(X)$ we denote the partition
${\bf 1}_X=\big\{ \{1,2,\dots,n\} \big\}$. This is the maximal element in
the set $\NC(X)$.

Suppose that $X=X_1 \cup X_2$, where $X_1\cap X_2=\emptyset$.
If $\pi\in\NC(X_1)$, we define its complementary partition
\index{complementary partition}
$\pi^{c}\in\NC(X_2)$ to be the biggest noncrossing partition of
$X_2$ such that $\pi\cup\pi^{c}\in\NC(X)$.

\subsubsection {Operator--valued noncrossing cumulants}
Let $(\D\subseteq\A,E)$ be an operator--valued probability space. The
algebra $\A$ has a structure of a $\D$--bimodule, we can define therefore
\begin{align}
\A^{\otimes_{\D} 0} & =\D, \notag \\
\A^{\otimes_{\D} n} &
=\underbrace{\A\otimes_{\D}\A\otimes_{\D}\cdots\otimes_{\D}\A}_{n\ {\rm
times}} \qquad \mbox{ for } n\geq 1 \notag
\end{align}
and a $\D$--bimodule
$$\Gamma=\bigoplus_{n\geq 0} \A^{\otimes_{\D} n}=\D\oplus \A\oplus
(\A\otimes_{\D}\A)\oplus\cdots, $$
called the free Fock space. \index{free Fock space}
We will denote $\Omega=I\oplus 0 \oplus 0 \oplus
\cdots \in \Gamma$. For every vector $v=d \oplus b_{1,1} \oplus (b_{2,1}
\otimes b_{2,2}) \oplus \cdots\in\Gamma$ we define $\langle \Omega,
v\rangle=d\in\D$ (it is possible to furnish the whole Fock space with a
structure of a Hilbert module but in this article we will use only scalar
products of the above form).

Let a $\D$--bimodule map $\kk:\bigoplus_{n\geq 1} \A^{\otimes_{\D}
n}\rightarrow \D$ be given. The map $\rho$ defined below assigns to each element of
$\A$ an operator acting on the Fock space $\Gamma$:
\begin{align} \rho(a)[d]&
=ad+\kk(a)d \qquad \qquad \mbox{for } d\in\A^{\otimes_{\D} 0}=\D, \notag
\end{align}
\begin{multline*}
\rho(a)[b_1\otimes \cdots \otimes b_n]
=a\otimes b_1 \otimes \cdots\otimes b_n  \\
+ \sum_{0\leq i\leq n} [\kk(a\otimes b_1 \otimes \cdots \otimes
b_i) b_{i+1}] \otimes b_{i+2} \otimes \cdots \otimes b_n
\end{multline*}
for $b_1 \otimes \cdots \otimes b_n\in\A^{\otimes_{\D} n}$.
Of course $\rho:\A\rightarrow L(\Gamma)$ is a $\D$--bimodule map.

We define a conditional expectation on the algebra $L(\Gamma)$ of linear
operators on $\Gamma$ by
$$\tilde{E}(x)=\langle\Omega,x\Omega\rangle.$$
In the literature this function is often called the vacuum expectation.

One can show \cite{Sp1998} that we can always choose the function $\kk$ so that
\begin{equation}
E(a_1 a_2 \cdots a_n)=\tilde{E}[\rho(a_1) \rho(a_2) \cdots
\rho(a_n)]
\label{eq:uwiklanemomenty}
\end{equation}
for every $n\geq 1$ and $a_1,\dots,a_n\in\A$; such unique
$\kk$ is called the operator--valued noncrossing cumulant function.
\index{operator--valued noncrossing cumulant function}
Sometimes in order to avoid ambiguities we shall denote it by $\kk^\D$.

\subsubsection {Noncrossing scalar cumulants and moments}
Let $(\A,\phi)$ be a scalar probability space.
Let $\kk:\bigoplus_{n\geq 1} \A^{\otimes_{\C}
n}\rightarrow\C$ be the $\C$--valued noncrossing cumulant function as
defined above. It will also be called scalar noncrossing cumulant and
sometimes we will denote it by $\kk^\C$.

For any $a_1,\dots\,a_n\in\A$
and a noncrossing partition $$\pi=\bigl\{ \{i_{1,1},i_{1,2},\dots,i_{1,m_1}
\}, \{i_{2,1},i_{2,2},\dots,i_{2,m_2} \}, \dots,
\{i_{k,1},i_{k,2},\dots,i_{k,m_k} \} \bigr\}$$
 of the set $\{1,2,\dots,n\}$ (where
$i_{j,1}<i_{j,2}<\cdots<i_{j,m_j}$ for every value of $j$) we define
\begin{equation}
\kk_{\pi}( a_1\otimes\cdots\otimes a_n)=
\prod_{1\leq l\leq k} \kk(a_{i_{l,1}}\otimes
a_{i_{l,2}}\otimes\cdots\otimes a_{i_{l,m_l}}).
\label{eq:kumulantyskalarne} \end{equation}
Cumulants $\kk$ will also be called simple cumulants;
\index{simple cumulants}
we will call $\kk_{\pi}$ compound cumulants.
\index{scalar compound cumulants}

It turns out that the implicit relation (\ref{eq:uwiklanemomenty}) between
cumulants and moments in the scalar case $\D=\C$ takes the form
\begin{equation} \phi(a_1
a_2\cdots a_n)=\sum_{\pi\in{\rm NC}(\{1,\dots,n\})} \kk_{\pi}
(a_1\otimes\cdots\otimes a_n)
\label{eq:momentyprzezkumulantyskalarne}
\end{equation}
for every $a_1,a_2,\dots,a_n\in\A$. In fact we can take this formula
as an alternative, inductive definition of $\kk$.

Similarly to (\ref{eq:kumulantyskalarne}) we define simple and compound
moments by
$$m(a_1\otimes \cdots \otimes a_n)=\phi(a_1 a_2\cdots a_n),$$
$$m_{\pi} (a_1\otimes\cdots\otimes
a_n)=\prod_{1\leq l\leq k} m(a_{i_{l,1}}\otimes a_{i_{l,2}}\otimes
\cdots\otimes a_{i_{l,m_l}}).$$

\begin{lemma} \label{lem:cyklicznosc}
Scalar cumulants have the following cyclic property
$$\kk^{\C}(a_1\otimes \cdots \otimes a_n)=
\kk^{\C} (a_n\otimes a_1 \otimes \cdots \otimes a_{n-1}).$$
\end{lemma}
The proof follows easily by induction from the tracial property of $\phi$.

\subsubsection {Noncrossing operator--valued cumulants and moments}
Let us consider an operator--valued probability space $(\D\subset\A,E)$.
The implicit relation (\ref{eq:uwiklanemomenty}) can be used to
express moments by operator--valued cumulants. For example we have:
\begin{align} \label{eq:przyklad0} E(a_1) & = \kk(a_1),\\
\label{eq:przyklad1}
E(a_1 a_2) & =\kk(a_1 \otimes a_2)+\kk(a_1)\kk(a_2), \end{align}
\begin{multline}
\label{eq:przyklad2} E(a_1 a_2 a_3)=
\kk(a_1 \otimes a_2 \otimes a_3)+ \kk(a_1 \otimes a_2)
\kk(a_3) \\
+\kk\big(a_1\otimes \kk(a_2) a_3 \big)+
\kk(a_1) \kk(a_2 \otimes a_3)+\kk(a_1) \kk(a_2) \kk(a_3).
\end{multline}
Summands on the right hand side correspond to all possible ways of writing
brackets in the expression $a_1 a_2 \cdots a_n$. In the general case the
following analogue of (\ref{eq:momentyprzezkumulantyskalarne}) holds
\begin{equation} E(a_1\cdots a_n)= \sum_{\pi\in{\rm NC}(\{1,\dots,n\})}
\kk_{\pi}( a_1\otimes \cdots \otimes a_n),
\label{eq:momentyprzezkumulanty} \end{equation} where the sum is taken over
all noncrossing partitions $\pi$ of the set $\{1,\dots,n\}$.
Unfortunately, since the algebra $\D$ does not commute with $\A$, the product
$\kk(a_{i_{1,1}}\otimes\cdots\otimes a_{i_{1,m_1}}) \cdots
\kk(a_{i_{k,1}}\otimes\cdots\otimes a_{i_{k,m_k}})$ in the definition of
compound cumulants (\ref{eq:kumulantyskalarne}) has to be replaced by a
nested product given by the partition $\pi$. For example we have (compare
with (\ref{eq:przyklad2})):
\begin{align*}
\kk_{\{\{1,2,3\}\}} (a_1\otimes a_2 \otimes a_3 ) & =
\kk (a_1\otimes a_2 \otimes a_3), \\
\kk_{\{\{1,2\},\{3\} \}} (a_1\otimes a_2 \otimes a_3) & =
 \kk (a_1 \otimes a_2) \kk(a_3), \\
\kk_{\{\{1,3\},\{2\} \}} (a_1\otimes a_2 \otimes a_3) & =
\kk\big(a_1\otimes \kk(a_2) a_3\big), \\
\kk_{\{ \{1\},\{2\},\{3\} \} }(a_1\otimes a_2\otimes a_3) & =
\kk(a_1) \kk(a_2) \kk(a_3).
\end{align*}

For a formal definition of compound operator--valued cumulants
\index{compound operator--valued cumulants}
$\kk_{\pi} (a_1\otimes\cdots\otimes a_n)$ we
refer to \cite{Sp1998}.

We can consider (\ref{eq:momentyprzezkumulanty}) as an alternative, inductive
definition of $\kk$.

Compound moments $E_{\pi}$ \index{compound moments}
are defined by the same recursive formula as
compound cumulants; for example
\begin{align*}
E_{\{\{1,2,3\}\}} (a_1\otimes a_2 \otimes a_3 ) & = E (a_1  a_2 a_3),\\
E_{\{\{1,2\},\{3\} \}} (a_1\otimes a_2 \otimes a_3)& =
 E (a_1  a_2) E(a_3),\\
E_{\{\{1,3\},\{2\} \}} (a_1\otimes a_2 \otimes a_3) &
= E\big(a_1 E(a_2) a_3\big), \\
E_{\{ \{1\},\{2\},\{3\} \} }(a_1\otimes a_2\otimes a_3) & =
E(a_1) E(a_2) E(a_3).
\end{align*}

\subsubsection {Freeness and cumulants}
The following theorem has been proved by Speicher \cite{Sp1998}.
\begin{theorem}
\label{theo:kumulantytowolnosc}
Let $(\D\subseteq\A,E)$ be an operator--valued probability space.
A family $(X_s)_{s\in S}$ of subsets of $\A$ is free with amalgamation
over $\D$ iff for every $n\geq 2$, $s_1,\dots,s_n\in S$ which are not all
equal, every $x_1\in X_{s_1},\dots,x_n\in X_{s_n}$, and every $d_1,\dots,d_n\in\D$
we have
$$\kk (x_1 d_1\otimes x_2 d_2\otimes \cdots \otimes x_{n} d_{n})=0.$$
\end{theorem}

Let us also recall the following characterization of the product
of free random variables, cf.\ \cite{NiSp1995}.
\begin{lemma} \label{lem:lematokumulantowaniu}
If $\{a_1,a_2,\dots,a_n,a_{n+1}\}\subseteq(\A,\phi)$ and
$\{b_1,b_2,\dots,b_n\}\subseteq(\A,\phi)$ are free then
\begin{multline}
\phi ([a_1] b_1 a_2 b_2\cdots a_n b_n
[a_{n+1}])\label{eq:momentyikumulanty} \\
=\sum_\pi
m_{\pi^{c}} ([a_1\otimes] a_2\otimes\cdots \otimes a_n[\otimes a_{n+1}])
\kk_{\pi} (b_1\otimes b_2\otimes \cdots \otimes b_n),
\end{multline}
 (the factors in brackets might be absent) where the sum is taken over all
noncrossing partitions $\pi$ of the $n$--element set of factors
$\{b_1,b_2,\dots,b_n\}$ and where $\pi^{c}$ denotes the complementary
partition of the set $\{[a_1,]a_2,a_3,\dots,a_n[,a_{n+1}]\}$.

Conversely, if at least one of two sets $X, Y\subseteq \A$ is a unital
subalgebra of $\A$ and (\ref{eq:momentyikumulanty}) holds for all
$a_1,\dots,a_{n+1}\in X$ and $b_1,\dots,b_n\in Y$ then $X$ and $Y$ are free.
\end{lemma}
One can write a $\D$--valued version the above lemma where freeness is
replaced by freeness with amalgamation and scalar cumulants are replaced by
operator--valued cumulants; the only thing we need to change is
to replace on the right--hand side the usual product by the nested
product given by the noncrossing partition $\pi^c\cup\pi$ of the
set of factors $\{[a_1,]b_1,a_2,b_2,\dots,a_n,b_n[,a_{n+1}]\}$.
For example we have
$$E[a_1  b_1  a_2 b_2] =E(a_1) \kk^\D \big(b_1\otimes E(a_2) b_2\big)+
E\big(a_1 \kk^\D(b_1) a_2\big) \kk^\D(b_2).$$
In general we have
\begin{multline} E([a_1] b_1 a_2
b_2 \cdots a_n b_n [a_{n+1}])= \\
\label{eq:momentyikumulanty2}
\sum_{\pi} (E\cup \kk^\D)_{\pi^{c}\cup\pi}([a_1\otimes]b_1\otimes
a_2\otimes\cdots\otimes a_n \otimes b_n [\otimes a_{n+1}]). \end{multline}
 We refer to \cite{Sp1998} for the concrete definition of $(E\cup
\kk^\D)_{\pi^c\cup\pi}$.

\subsubsection {Theorem of Krawczyk and Speicher}
We will use often the following theorem of Krawczyk and Speicher
\cite{KS} which allows us to compute cumulants of products of random
variables.

\begin{theorem} \label{theo:krawczykspeicher}
Let $m\in\N$ and $1\leq i_1<i_2<\cdots <i_m:=n$ be given. Consider random
variables $a_1,\dots,a_n$ and put $A_j:=a_{i_{j-1}+1} \cdots a_{i_j}$
for $j=1,\dots,m$ (where $i_0:=0$).
Then the following equation holds:
\begin{equation} \label{eq:krawczykspeicher}
\kk [(a_1 \cdots a_{i_1})\otimes\cdots\otimes
( a_{i_{m-1}+1} \cdots a_{i_m})]=
\sum_{
{{\pi\in\NC(\{1,\dots,n\})} \atop {\pi\vee\sigma={\mathbf 1}_n}}
}
\kk_{\pi} [a_1\otimes\cdots\otimes a_n],
\end{equation}
where $\sigma\in \NC(\{1,\dots,n\})$ is the partition
$\sigma=\big\{ \{1,\dots,i_1\},\dots,\{i_{m-1}+1,\dots,i_m\} \big\}$.
\end{theorem}

This theorem was originally proved in the scalar case, but the proof
in \cite{KS} extends without difficulties to the operator--valued setup.

\subsubsection {Scalar random variables. Scalar Haar unitaries.}
We say that $y\in(\D\subseteq\A,E)$ is a scalar random variable
\index{scalar random variable} if $yy\gwia=y\gwia y$,
$y$ and $y\gwia$ commute with elements of $\D$ and if for every
$k,l\in\{0,1,2,\dots\}$ we have
$E(y^k (y\gwia)^l)=\alpha_{k,l} I$ for some $\alpha_{k,l}\in\C$.

We say that $u\in(\D\subseteq\A,E)$ is a scalar Haar unitary
\index{scalar Haar unitary} if $u$ is a
scalar random variable which is a unitary such that $E(u^n)=E[(u\gwia)^n]=0$
for every $n\in\N$, $n\geq 1$.

\subsection {Voiculescu's $R$-- and $S$--transforms}
In this section we recall basic properties of Voiculescu's
$R$-- and $S$--transforms
(cf \cite{Voi1986,Voi1987,NiSp1994,Nica1996,Ha1997,HP}).

In the following $(\cdot)^{\langle -1\rangle}$ will denote
the inverse function with respect to composition.

For an element $a$ of a scalar probability space $(\A,\phi)$ we
define its Cauchy transform \index{Cauchy transform} $G_a$ by
$$G_a(t)=\phi\left( \frac{1}{t-a} \right) \qquad\mbox{for }
t\not\in\sigma(a).$$
The function
$$R_a(z)=G^{\langle -1\rangle}(z)-\frac{1}{z} $$
makes sense in a neighborhood of $0$, is an analytic function
and is called the $R$--transform of $a$. \index{$R$--transform}

We define furthermore
$$\psi_a(t)=\phi\left(\frac{1}{1-ta}\right)-1 \qquad \mbox{for }
\frac{1}{t}\not\in\sigma(a).$$
If $\phi(a)\neq 0$ then
$$S_a(z)=\frac{z+1}{z} \psi^{\langle -1\rangle}(z)$$
makes sense in a neighborhood of $0$, is an analytic function
and is called the $S$--transform of $a$. \index{$S$--transform}

$R$--diagonal elements will be defined in the next Section , despite this we
will state the following lemma here.

\begin{lemma} \label{lem:stransformataipromienspektralny} Let $x$ and $y$
($x,y\neq 0$) be $R$--diagonal elements.
The spectral radius $y$ is given by
$$\sup_{\lambda\in\sigma(y)} |\lambda|=\|y\|_{L^2}=S_{(yy\gwia)^{-1}}(-1).$$

The spectral radius
of $x^{-1}$ is given by $S_{x x\gwia}(-1)$, or equivalently,
$$\inf_{\lambda\in\sigma(x)} |\lambda|=\frac{1}{S_{x x\gwia}(-1)}.$$
\end{lemma}
\begin{proof}
The identity $\sup_{\lambda\in\sigma(y)} |\lambda|=\|y\|_{L^2}$
was proved by Haagerup and Larsen \cite{HL,La}.

We will prove now that $S_a(-1)=\phi(a^{-1})$ holds for every
nonnegative operator $a$. Our method is based on proofs of the formula
\begin{equation} z S_a(z)=[z
R_a(z)]^{\langle -1 \rangle}. \label{eq:dziwnaformulanaRiS}
\end{equation}
from articles \cite{NiSp1997A,HL}.

If $a$ is a nonnegative operator, then the function
$\psi_a:(-\infty,0]\rightarrow\R$ is a well--defined increasing function
with the range $\big( \phi(p_0)-1,0 \big]$, where $p_0$ denotes the
projection on the kernel of $a$. It means that the $S$--transform of $a$
is a well defined analytical function on $\big( \phi(p_0)-1,0 \big)$.

If $\phi(p_0)>0$ then $S_a(-1)$ is not well defined and we define
$S_a(-1)=\infty$, therefore $S_a(-1)=\phi(a^{-1})$ holds by definition.

Let us consider now the case $\phi(p_0)=0$; we take
$$S_a(-1):=\lim_{\epsilon\rightarrow 0^{+}} S_a(-1+\epsilon)$$
if this limit exists and $S_a(-1)=\infty$ otherwise.

For every $u$ such that $u^{-1}\not\in\sigma(a)$ we have that
$$\psi_a (u)=\frac{G_a(u^{-1})}{u}-1,$$
so for every $u<0$ we have
$$\psi_a(u) S_a[\psi_a(u)]=[1+\psi_a(u)] \psi^{\langle -1\rangle}_a[\psi_a(u)]=
 [1+\psi_a(u)] u =G_a(u^{-1}).$$
Take the limit $u\rightarrow -\infty$; since $\phi(p_0)=0$ we have
$\lim_{u\rightarrow -\infty} \psi_a(u)=-1$ and therefore
$$S_a(-1)=\phi(a^{-1}).$$

%
\end{proof}

\begin{lemma} \label{lem:stransformatakompresji}
Let $R^{\A}_a$ and $S^{\A}_a$ denote the $R$-- and $S$--transform, respectively,
of an element $a$ in the probability space $(\A,\phi)$. Let $q$ be an
orthogonal projection which is free from $a$ and which fulfills
$0<\phi(q)=s$. Let $R^{q \A q}_{q a q}$ and $S^{q \A q}_{q a
q}$ denote $R$-- and $S$--transform, respectively, of the element
$q a q\in(q \A q, \frac{1}{s} \phi)$. The
following identities hold:
$$R^{q \A q}_{q a q}(t)=R^{\A}_{a}(st), $$
$$S^{q \A q}_{q a q}(z)=S^{\A}_{a}(sz).$$
\end{lemma}
\begin{proof}
The first identity was proved in papers \cite{NiSp1995,Dima1997}.
The second identity follows from the first one and (\ref{eq:dziwnaformulanaRiS}).
\end{proof}

The following lemma is due to Haagerup and Larsen \cite{HL}.
The proof of the first identity is by direct calculation and the second identity follows
from the first one.
\begin{lemma} \label{lem:stransformataodwrotnego}
For every noncommutative random variable $a$ and $t\not\in\sigma(a)$
we have
$$\psi_{a^{-1}}(t)+\psi_a\left(\frac{1}{t}\right)=-1.$$

If $a$ is a strictly positive operator and $0\leq s\leq 1$ then
$$S_a(-s) S_{a^{-1}}(s-1)=1.$$
\end{lemma}


\section {$R$--diagonality with amalgamation}
\label{chap:rdiagonality}
\subsection {Definition and basic properties}
\label{sec:rdiagonality}
In this section we show a few equivalent definitions of $R$--diagonality
with amalgamation. Usually they are direct analogues of their scalar
analogues; proof of their equivalence is similar to the one in the scalar
case as well. These definitions are
expressed in terms of different notions, what will turn out to be very
useful in applications.

\begin{thedef} \label{thedef:centralne}
Let $(\D\subseteq\A,E)$ be a $\D$--valued probability space and let
$a\in\A$. Then the following conditions on $a$ are equivalent.
The element $a$ is said to be $R$-diagonal with amalgamation over $\D$
(or $R$--diagonal over $\D$)
\index{$R$--diagonality with amalgamation} if it
satisfies one (hence all) of these conditions.
\begin{enumerate}

\item \label{war:momenty}
We define
\begin{align*}
P_{11} & = \spann\{ \D a\gwia \D, \D a\gwia \D a \D a\gwia \D, \D a\gwia \D a
\D a\gwia \D a \D a\gwia \D ,\dots\},\\
P_{12} & =\big\{ x\in\spann \{ \D, \D a\gwia\D a \D,
\D a\gwia \D a \D a\gwia \D a \D, \dots \} : E(x)=0 \big\}, \\
P_{21} & =\big\{ x\in\spann \{\D, \D a \D a\gwia \D, \D a \D a\gwia \D a \D
a\gwia \D,\dots \}: E(x)=0 \big\}, \\
P_{22} & =\spann \{ \D a \D, \D a \D a\gwia \D a \D, \D a \D a\gwia
\D a \D a\gwia\D a\D, \dots \}.
\end{align*}

One has
$$E(x_1 x_2 \cdots x_n)=0,$$
for every $n\geq 1$, all $i_1,i_2,\dots,i_{n+1}\in\{1,2\}$, and all $x_1\in
P_{i_1 i_2}$, $x_2\in P_{i_2 i_3}$, \dots, $x_{n}\in P_{i_n,i_{n+1}}$.

\item \label{war:kumulanty}
Consider the family of $\D$--valued noncrossing cumulants
$$\bigl\{ \kk(d_1 a^{s_1}\otimes d_2 a^{s_2}\otimes\cdots\otimes d_n
a^{s_n}): n\geq 1; s_1,\dots,s_n\in\{1,\star\}, d_1,\dots,d_n\in\D
\bigr\}.$$ Then $ \kk(d_1 a^{s_1}\otimes d_2 a^{s_2}\otimes\cdots\otimes
d_n a^{s_n})=0$ whenever $(s_1,\dots,s_n)$ is not alternating, i.e.\ if
it is not of the form
$(1,\star,1,\star,\dots,1,\star)$ or $(\star,1,\star,1,\dots,\star,1)$.

\item \label{war:istniejeunitarny}
There exists an enlargement $(\D\subseteq\tilde{\A},\tilde{E})$ of the
$\D$--valued probability space $(\D\subseteq\A,E)$ and a scalar element
$u\in\tilde{\A}$ such that:
\begin{enumerate} \item $u$ is unitary;
\item $\{a,a\gwia\}$ and
$\{u,u\gwia\}$ are free with amalgamation over $\D$;
\item $|E(u)|<1$ (since $u$ is
scalar, $E(u)$ is a multiple of identity, can be therefore regarded as
a complex number);
\item elements $a$ and $ua$ are identically $\D$--distributed.
\end{enumerate}

\item \label{war:kazdyunitarny}
For every enlargement
$(\D\subseteq\tilde{\A},\tilde{E})$ of the $\D$--valued probability space
$(\D\subseteq\A,E)$ and a scalar element $u\in\tilde{\A}$ such that $u$ is
unitary and $\{a,a\gwia\}$ and $\{u,u\gwia\}$ are free with amalgamation
over $\D$, we have that elements $a$ and $ua$ are identically
$\D$--distributed.

\end{enumerate}

If $x\in(\A,\phi)$ is $R$--diagonal with amalgamation over $\C$ we simply say
that $x$ is $R$--diagonal (or scalar $R$--diagonal). \index{$R$--diagonality}
\end{thedef}

\begin{proof}
The proof follows closely the proof the corresponding scalar--valued
Theorem from \cite{NiShSp} and thus we will be quite condensed in the
following.

{\bf Proof of the implications $\ref{war:istniejeunitarny}\Longrightarrow
\ref{war:momenty}\Longrightarrow \ref{war:kazdyunitarny}$}

Suppose that $u$ is a scalar unitary. Let $x\mapsto x'$ be a linear
map which acts on formal noncommutative polynomials of $a$ and $a\gwia$ with
coefficients in $\D$ and which replaces $a$ by $ua$:
$$d_1
a^{s_1} d_2 a^{s_2} \cdots d_{n} a^{s_n} d_{n+1} \mapsto d_1 (ua)^{s_1} d_2
(ua)^{s_2} \cdots d_{n} (ua)^{s_n} d_{n+1},$$ for $n\geq 0$,
$d_1,\dots,d_{n+1}\in\D$, $s_1,\dots,s_n\in\{1,\star\}$.

It is easy to see that every polynomial of $a$, $a\gwia$ and $\D$ is
a linear combination of an element $d\in\D$ and expressions of type
\begin{multline} \label{eq:wyrazenie}
x=x_1 x_2 \cdots x_n, \quad
\mbox{where } n\geq 1; \\ x_1\in P_{i_1 i_2},
x_2\in P_{i_2 i_3}, \dots ,x_n\in P_{i_n i_{n+1}};\
 i_1,\dots, i_{n+1}\in\{1,2\},
\end{multline}
therefore the
property that $a$ and $ua$ are identically $\D$--distributed is
equivalent to the property that for every $x$ as in Eq.\
(\ref{eq:wyrazenie}) we have $E(x')=E(x)$.

We have that
\begin{align*}
 y' &=y u\gwia           && \quad \mbox{for } y\in P_{11}, \\
 y' &=y                  && \quad \mbox{for } y\in P_{12}, \\
 y' &=uyu\gwia           && \quad \mbox{for } y\in P_{21}, \\
 y' &=u y                &&  \quad \mbox{for } y\in P_{22},
\end{align*}
therefore for $x$ as in (\ref{eq:wyrazenie})
there are
$s_1,\dots,s_{n+1}\in\{1,\star\}$ such that
\begin{equation} \label{eq:wzornaxprim}
x'=x_1' x_2' \cdots
x_n'=[u^{s_1}] x_1 u^{s_2} x_2 \cdots u^{s_{n}} x_n [u^{s_{n+1}}]
\end{equation} 
(the factors in brackets might be absent).

($\ref{war:momenty}\Longrightarrow \ref{war:kazdyunitarny}$)
In Eq.\ (\ref{eq:wzornaxprim}) we can replace $u$ by
$E(u)+u_0$, where $u_0=u-E(u)$, and distribute it. Since $E(u_0)=0$ and
$E(u)\in\C$, for each summand we can directly apply the freeness condition
(\ref{eq:freeness}), therefore $E(x')=0=E(x)$ and condition
\ref{war:kazdyunitarny} holds.

($\ref{war:istniejeunitarny}\Longrightarrow \ref{war:momenty}$)
Let $u_1,u_2,\dots$ be a sequence of scalar unitaries identically distributed
as $u$, such that the family $\{a,a\gwia\},\{u_1,u_1\gwia\},\{u_2,u_2\gwia\},
\dots$ is free with amalgamation over $\D$
(it is always possible to extend the probability space
$(\D\subseteq\tilde{\A},\tilde{E})$ in order to find such unitaries).
It follows by induction
that $a$ and $v_n a$ have the same $\D$--valued distributions, where
$v_n=u_n u_{n-1} \cdots u_2 u_1$ are scalar unitaries. For each $n$ we have
that $\{v_n,v_n\gwia\}$ and $\{a,a\gwia\}$ are free and $E(v_n)$ converges
to zero. It follows that without any loss of generality we can assume that
$E(u)=0$.

We have that if $z\in P_{11}$ then $E(z)=E(z')=E(z u)=E(z) E(u)$ and
therefore $E(z)=0$. Similarly if $z\in P_{22}$ then $E(z)=0$; therefore if
$z\in P_{ij}$, $i,j\in\{1,2\}$ then $E(z)=0$.

{} From (\ref{eq:wzornaxprim}) and the freeness condition (\ref{eq:freeness})
we have $$E(x)=E(x')=E([u^{s_1}] x_1 u^{s_2} x_2 \cdots u^{s_{n}} x_n
[u^{s_{n+1}}])=0,$$
therefore condition $\ref{war:momenty}$ holds.

{\bf Proof of the equivalence $\ref{war:kumulanty}\iff
\ref{war:kazdyunitarny}$}

($\ref{war:kumulanty}\Longrightarrow \ref{war:kazdyunitarny}$) Let
us consider one of the moments
\begin{equation} E[(a u)^{s_1}d_1 (a u)^{s_2} d_2 \cdots d_{n-1}
(a u)^{s_n}] , \label{eq:moment} \end{equation} where
$s_1,\dots,s_n\in\{1,\star\}$ and $d_1,\dots,d_{n-1}\in\D$. It can
be written in the form (\ref{eq:momentyikumulanty2}), where
$a_1,\dots,a_{n+1}\in \{u,d,u\gwia: d\in\D \}$ and
$b_1,\dots,b_n\in\{d a, a\gwia d: d\in\D\}$. Simple combinatorial
arguments show that only special non--crossing partitions
contribute in Eq.\ (\ref{eq:momentyikumulanty2}), namely
noncrossing partitions $\pi$ such that $E_{\pi^{c}}$
does not depend on $u$ and condition $\ref{war:kazdyunitarny}$ holds.

($\ref{war:kazdyunitarny}\Longrightarrow  \ref{war:kumulanty}$)
Suppose that the statement holds for every $n<n_0$ and let us consider
$(s_1,s_2,\dots,s_n)$, $s_1,\dots,s_n\in\{1,\star\}$, such that
$(s_1,s_2,\dots,s_n)$ is equal neither to $(1,\star,1,\star,\dots,1,\star)$
nor to $(\star,1,\star,1,\dots,\star,1)$.

As above, we can write (\ref{eq:moment})
in the form (\ref{eq:momentyikumulanty2}), where
$a_1,\dots,a_{n+1}\in \{u,d,u\gwia: d\in\D \}$ and
$b_1,\dots,b_n\in\{d a, a\gwia d: d\in\D\}$. From the inductive assumption
we see that the only summand on the right hand side of
(\ref{eq:momentyikumulanty2}) which depends on $E(u)$ or $E(u\gwia)$ is
the summand corresponding to $\pi=\big\{\{1,2,\dots,n\}\big\}$; therefore
it must be equal to zero and the statement follows.

Of course the implication $\ref{war:kazdyunitarny}
\Longrightarrow \ref{war:istniejeunitarny}$ is  trivial.
\end{proof}

\begin{remark}
\label{newremark}
In the same way as in the scalar--valued case \cite{NiSp1997A}, it is easy
to see the following:
\begin{enumerate}
\item Scalar Haar unitaries are $R$--diagonal with amalgamation over $\D$.
\item If $a,b\in(\D\subseteq\A,E)$ have a property that
 $\{a,a\gwia\}$ and $\{b,b\gwia\}$ are free with amalgamation
over $\D$ and $a$ is $R$--diagonal with amalgamation over $\D$
then also $ab$ and $ba$ are $R$--diagonal with amalgamation over
$\D$.
\end{enumerate}
\end{remark}

\subsection {Scalar $R$--diagonality in terms of freeness}
We say that a nonempty
sequence $(s_1,\dots,s_n)$ (where $s_1,\dots,s_n\in\{-1,1\}$)
is balanced \index{balanced sequence}
if $s_1+\cdots+s_n=0$ and its partial sums fulfill
$s_1,s_1+s_2,\dots, s_1+\cdots+s_{n-1}\leq 0$.

Let an element $a\in(\A,\phi)$ be given. We will use the following notation:
$a^1:=a$ and $a^{-1}:=a\gwia$.
We define \index{algebra $K_a$}
$$\K_a=\spann\{ a^{s_1} a^{s_2} \cdots a^{s_n}: (s_1,\dots,s_n)
\mbox{ is balanced}\}.$$

It is easy to see that $\K_a$ is an algebra.

\begin{theorem} \label{theo:nowardiagonalnosc}
An element $a\in(\A,\phi)$ is
$R$--diagonal if and only if the following conditions hold:
\begin{enumerate}
\item \label{war:1} for every $z\in\C$ such that $|z|=1$ we have
that $a$ and $za$ have the same distribution,
\item \label{war:2} $a a\gwia$ and $K_a$ are free.
\end{enumerate}
\end{theorem}
\begin{proof}
Suppose that $a$ is $R$--diagonal. It trivially fulfills condition
\ref{war:1}.

We use Theorem \ref{theo:krawczykspeicher} to compute a cumulant
$\kk(x_1 \otimes \cdots \otimes x_n)$, where for each $1\leq k\leq n$ we have that
either $x_k=aa\gwia$ or $x_k=a^{s^k_1} \cdots a^{s^k_{n_k}}$,
where the sequence $(s^k_i)_{1\leq i\leq n_k}$ is balanced; furthermore
let each of these two possibilities occur for at least one index $k$.
From the following Lemma \ref{lem:lematogorkowaniu} it follows that
the sum on the right--hand side of (\ref{eq:krawczykspeicher}) is equal to
zero; therefore
Theorem \ref{theo:kumulantytowolnosc} implies that the
condition \ref{war:2} is fulfilled.

Suppose now that $a\in\A$ fulfills conditions \ref{war:1} and \ref{war:2}.
We will prove by induction over $n$ that
$\kk(a^{s_1}\otimes\cdots\otimes a^{s_n})=0$ for
all sequences $s_1,\dots,s_n\in\{-1,1\}$ which are not alternating, i.e.\
are not of the following form: $(-1,1,-1,1,\dots,-1,1)$ or
$(1,-1,1,-1,\dots,1,-1)$.

Consider first the case $s_1+\cdots+s_n\neq 0$. From condition \ref{war:1}
follows that for every $z\in\C$, $|z|=1$ we have
\begin{multline} \kk(a^{s_1}\otimes \cdots \otimes a^{s_n})=
\kk(z^{s_1} a^{s_1} \otimes z^{s_2} a^{s_2} \otimes \cdots
\otimes z^{s_n} a^{s_n})\\
=z^{s_1+\cdots+s_n} \kk(a^{s_1}\otimes \cdots \otimes a^{s_n});
\end{multline}
therefore  $\kk(a^{s_1}\otimes \cdots \otimes a^{s_n})=0$.

We consider now the case $s_1+\cdots+s_n=0$. We define $k$ to be the index
such that the partial sum $s_1+\dots+s_{k+1}$ takes the biggest value and
a sequence $(r_i)$ given by a cyclic rotation of the sequence $(s_i)$:
$$(r_1,r_2,\dots,r_n)=
(s_{n-k+1},s_{n-k+2},\dots,s_{n-1},s_n,s_1,s_2,\dots,s_{n-k}).$$
Since the scalar--valued cumulants $\kk$ have the cyclic property
(cf Lemma \ref{lem:cyklicznosc})
it is sufficient to show our theorem for the sequence $(r_i)$.

Partial sums of $(r_i)$ have the following property:
$$r_1,r_1+r_2, r_1+r_2+r_3, \dots, r_1+\cdots+r_{n-1} \leq 1$$
and $r_1=1$. It follows that $(r_i)$ can be written as a concatenation
of the following nonempty sequences:
$$(r_1,\dots,r_n)=\big(  (s^1),  (s^2), \dots,  (s^m) \big),$$
where for each $k$ we have that either the sequence $(s^k)$ is given by
$(s^k)=(1,-1)$ or $(s^k)=(s^k_1, s^k_2, \dots, s^k_{n_k})$ is balanced.
Every possibility occurs at least once unless the original sequence
$(s_1,\dots,s_n)$ was not equal to $(1,-1,1,-1,\dots,1,-1)$ or to
$(-1,1,-1,1,\dots,-1,1)$.
Therefore the condition \ref{war:2} implies that
\begin{equation} \kk\big( (a^{s^1_1} \cdots a^{s^1_{n_1}})
\otimes \cdots \otimes (a^{s^m_1}\cdots a^{s^m_{n_m}}) \big)=0.
\label{eq:1410}
\end{equation}

On the other hand we can use Theorem \ref{theo:krawczykspeicher} to
express the left--hand side of (\ref{eq:1410}) in terms of cumulants of type
$\kk(a^{p_1}\otimes \cdots \otimes a^{p_w})$.
At the first sight Lemma \ref{lem:lematogorkowaniu} concerns only
$R$--diagonal operators; however the inductive hypothesis assures us that
all non--alternating cumulants of less than $n$ factors vanish;
therefore we have that the only nonvanishing term in
(\ref{eq:krawczykspeicher}) is given by the full partition $\pi=\jed_n$:
\begin{equation} \kk\big( (a^{s^1_1} \cdots a^{s^1_{n_1}})
\otimes \cdots \otimes (a^{s^m_1}\cdots a^{s^m_{n_m}}) \big)=
\kk( a^{r_1} \otimes \cdots \otimes a^{r_n}).
\label{eq:1411}
\end{equation}

Eq.\ (\ref{eq:1410}) and (\ref{eq:1411}) imply that
$\kk( a^{r_1} \otimes \cdots \otimes a^{r_n})=0$.
\end{proof}


The above proof used the tracial property of the state $\phi$, therefore
it seems to be technically difficult to find a similar characterization
of the $R$--diagonality in the operator--valued setup.
However, we have the following implication in one direction:
\begin{theorem} \label{theo:nowardiagonalnoscoperatorowa}
Suppose that $a\in(\D\subset\A,E)$ is $R$--diagonal with amalgamation.
Then the following hold:
\begin{enumerate}
\item for every $z\in\C$ such that $|z|=1$ we have that $a$ and $za$ have
the same distribution,
\item $a a\gwia$ and $K_a$ are free with amalgamation over $\D$.
\end{enumerate}
\end{theorem}
The proof follows exactly its scalar--valued analogue.

\begin{lemma} \label{lem:lematogorkowaniu}
Suppose that $a\in(\A,\phi)$ is $R$--diagonal, let
$t_1,\dots,t_n\in\{-1,1\}$. We introduce an equivalence relation on the
set $\{1,\dots,n\}$ given by
$$(i \sim j) \iff t_1+\cdots+t_{i-1}+\frac{t_i}{2}=t_1+\cdots+t_{j-1}+\frac{t_j}{2}.$$
The above equivalence relation defines a certain (possibly crossing) partition
$\sigma$ of the set $\{1,\dots,n\}$.

Let $\pi\in\NC(\{1,\dots,n\})$ be such that
$\kk_{\pi}(a^{t_1} \otimes \cdots \otimes a^{t_n})\neq 0$.

Then $t_1+\cdots+t_n=0$ and $\pi\prec\sigma$.
\end{lemma}
The proof follows easily by induction.

\section {Compatible scalar and operator--valued structures}
\label{chap:compatible}
In the following we will consider an algebra $\A$ equipped with two structures:
the scalar--valued structure $(\A,\phi)$ and an operator--valued structure
$(\D\subseteq\A,E)$. We say that the scalar and operator--valued structures on
$\A$ are compatible if $\phi(x)=\phi[E(x)]$ for every $x\in\A$.
\index{compatible structures}.

It is a fascinating question how to relate properties of random variables
regarded once as elements of $(\A,\phi)$ and another tome as elements of
$(\D\subseteq\A,E)$. Of course the $\D$--valued distribution (with respect to
conditional expectation $E$) of a random variable $x$ determines its
$\C$--valued distribution (with respect to the trace $\phi$), but this
formula might be very complicated. In the following we consider some special
situations, when it is much easier to establish such a link.

In order to distinguish these objects, we will denote by $\kk^{\C}$
the cumulants $\kk:\A^{\otimes n}\rightarrow \C$ in the scalar
probability space $(\A,\phi)$ and by $\kk^{\D}$
the cumulants $\kk:\A^{\otimes_{\D} n}\rightarrow \D$ in the $\D$--valued
probability space $(\D\subseteq\A,E)$.

\subsection {Freeness with algebra $\D$}
The following theorem is a special case of a result of Nica, Shlyakhtenko, and
Speicher, which will be published elsewhere
\cite{NiShSp3}. In order to be self-contained we
provide a short sketch of the proof for our special case.
\begin{theorem}
\label{theo:wolneodskalarow2}
Let $(\A,\phi)$ and $(\D\subseteq\A,E)$ be compatible and let $X\subseteq\A$.
Then, $X$ and $\D$ are free in $(\A,\phi)$ iff
for every $n\geq 1$ and $x_1,\dots,x_n\in X$ there exists
$c_n(x_1,\dots,x_n)\in\C$
such that for every $d_1,\dots,d_{n-1}\in\D$ we have
\begin{multline} \kk^{\D}(x_1 d_1 \otimes x_2 d_2\otimes\cdots\otimes
x_{n-1} d_{n-1}\otimes x_n)\\ =c_n(x_1,\dots,x_n)
\phi(d_1) \phi(d_2) \cdots \phi(d_{n-1}) I.
\label{eq:1978} \end{multline}

If the above holds, we have
$$c_n(x_1,\dots,x_n)=\kk^{\C}(x_1 \otimes \cdots \otimes x_n).$$
\end{theorem}

\begin{proof}
Suppose that $X$ and $\D$ are free with respect to the state $\phi$.
We will prove by induction over $n$ that
\begin{multline*} \kk^{\D}(x_1 d_1 \otimes x_2 d_2\otimes\cdots\otimes x_{n-1}
d_{n-1}\otimes x)\\ =\kk^{\C}(x_1\otimes\cdots\otimes x_n)
\phi(d_1) \phi(d_2) \cdots \phi(d_{n-1}) I.\end{multline*}

Suppose that the statement is true for all $n<n_0$. Lemma
\ref{lem:lematokumulantowaniu} gives us
\begin{equation} \phi(x_1 d_1 x_2 d_2 \cdots
x_{n_0} d_{n_0})= \sum_{\pi}
\kk^{\C}_{\pi}(x_1\otimes \cdots\otimes x_{n_0})
m_{\pi^{c}}(d_1\otimes \cdots \otimes d_{n_0}) \label{eq:1976a}
\end{equation}
Directly from the definition of the free cumulants we have
\begin{equation} \phi[E(x_1 d_1 \cdots x_{n_0} d_{n_0})]= \sum_{\pi}
\phi[\kk^{\D}_{\pi}(x_1 d_1\otimes \cdots \otimes x_{n_0}
d_{n_0})]. \label{eq:1976b} \end{equation}
The left--hand sides of these equations are equal.
{} From the inductive hypothesis we see that if
$\pi\neq\big\{\{1,2,\dots,n\}\big\}$ then corresponding summands on the
right--hand sides of (\ref{eq:1976a}) and (\ref{eq:1976b}) are
equal; therefore also for $\pi=\big\{\{1,2,\dots,n\}\big\}$
they must be equal as well:
\begin{multline*}
\phi\big[k^{\C}(x_1\otimes\cdots\otimes x_{n_0})
\phi(d_1)\cdots \phi(d_{n_0 -1})
d_{n_0}\big]=\\
\phi[\kk^{\D}(x_1 d_1\otimes \cdots
\otimes x_{n_0-1} d_{n_0-1} \otimes x_{n_0}) d_{n_0}]. \end{multline*}
Since $\phi$ is
faithful and the above identity holds for all $d_{n_0}\in\D$,
it follows $$k^{\C}(x_1\otimes\cdots\otimes x_{n_0})
\phi(d_1)\cdots \phi(d_{n_0-1})=\kk^{\D}(x_1
d_1\otimes \cdots \otimes x_{n_0-1} d_{n_0-1} \otimes x_{n_0}).$$

Conversely, if (\ref{eq:1978}) holds for all $n\in\N$, let us define
$\kk^{\C}(x_1\otimes\cdots\otimes x_n)=c_n(x_1,\dots,x_n)$. We have that
\begin{multline*}
\phi(x_1 d_1 x_2 d_2 \cdots x_{n} d_{n})=
\phi[E(x_1 d_1 x_2 d_2 \cdots x_{n} d_{n})]= \\
\sum_{\pi} \phi[\kk^{\D}_{\pi}(x_1 d_1\otimes \cdots \otimes
x_{n} d_{n})]=
 \sum_{\pi}
\kk^\C_{\pi}(x_1\otimes\cdots\otimes x_{n})
m_{\pi^{c}}(d_1\otimes \cdots \otimes d_{n})
\end{multline*}
If we take $d_1=\cdots=d_{n}=I$ we see that $\kk^{\C}$ indeed fulfills the
definition of a cumulant. Lemma \ref{lem:lematokumulantowaniu} implies
that $X$ and $\D$ are free.
%
%
\end{proof}

\subsection {Freeness versus freeness with amalgamation}
The following theorem is a special case of a result of Schlyakhtenko.
\begin{theorem}
\label{theo:schlyakhtenko}
Let $(\A,\phi)$ and $(\D\subseteq\A,E)$ be compatible,
let unital algebras $X$ and $\B$ (where $X,\B\subseteq\A$)
be free with amalgamation over $\D$
and let $X$ and $\D$ be free (with respect to
state $\phi$).

Then $X$ and $\D$ are free (with respect to the state $\phi$).
\end{theorem}
\begin{proof}
From the operator--valued version of Lemma \ref{lem:lematokumulantowaniu}
it follows that for $x_1,\dots,x_n\in X$ and $b_1,\dots,b_n\in\B$ we have
$$E[b_1 x_1 b_1 x_2 \cdots b_n x_n]=
\sum_{\pi} (E\cup \kk^{\D})_{\pi^{c}\cup \pi} (b_1\otimes x_1\otimes \cdots
\otimes b_n \otimes x_n).$$
We apply state $\phi$ to both sides of the equation; Theorem
\ref{theo:wolneodskalarow2} implies that
$$\phi(b_1 x_1\cdots b_n x_n)=\sum_{\pi} (\phi\cup \kk^\C)_{\pi\cup\pi^{c}}
(b_1\otimes x_1\otimes \cdots \otimes b_n \otimes x_n).$$

This means (Lemma \ref{lem:lematokumulantowaniu}) that the joint moments
of $X$ and $\B$ are as if $X$ and $\B$ were free. But this means exactly that
they are free.
\end{proof}

\subsection [Relation between scalar $R$-diagonality\dots]{Relation
between scalar $R$--diagonality and $R$--diagonality with
amalgamation}
\label{sec:relation}

The following theorem is one of the main results of this article.

\begin{theorem} \label{theo:rdiagonalnosc}
Let us consider an $\F$-valued probability space $(\F\subseteq \A,E)$
and a scalar probability space $(\A, \phi)$ which are compatible.
Suppose that $a\in\A$ is $R$-diagonal with amalgamation over $\F$ and
furthermore let $a a\gwia$ and $\F$ be free.
Then $a$ is $R$-diagonal.
\end{theorem}
We provide two proofs of this theorem.
\begin{proof}
We apply Theorem \ref{theo:nowardiagonalnoscoperatorowa} and
Theorem \ref{theo:schlyakhtenko} for algebra $X$ generated by $a a\gwia$.
It follows that the assumptions of Theorem \ref{theo:nowardiagonalnosc} are
fulfilled.
\end{proof}

\begin{proof}
Let $\B$ be the unital algebra generated by $aa\gwia$ and let $\B_0=\{x\in\B:
\phi(x)=0\}$. Furthermore we define $\F_0=\{x\in\F: \phi(x)=0\}$.

We define
\begin{multline*} Q_{11}={\rm span}\{(a\gwia \B), (a\gwia \B)
\F_0 \B_0, (a\gwia \B) \F_0 \B_0 \F_0 \B_0, \\
(a\gwia\B) \F_0 \B_0 \F_0 \B_0 \F_0 \B_0,\dots\},\end{multline*}
\begin{align*}
Q_{12} & =\bigl\{ x\in{\rm span} \{ \F,  a\gwia \B a ,
a\gwia \B \F \B a,a\gwia \B \F \B \F \B a, \dots \}: \phi(x)=0
\bigr\},\\
Q_{21} & ={\rm span} \{ \B_0,\B_0 \F_0 \B_0, \B_0 \F_0 \B_0
\F_0 \B_0,\dots \}, \\
Q_{22} & ={\rm span} \{ (\B a ), \B_0 \F_0 (\B a ), \B_0 \F_0
\B_0 \F_0 (\B a),\dots \}. \end{align*}

Let sets $P^{\F}_{ij}$ and $P^{\C}_{ij}$ be as in condition
\ref{war:momenty} of Theorem \ref{thedef:centralne} for $\D=\F$, resp.\ for
$\D=\C$ (in this case we use $\phi$ instead of $E$).

We have that $P^{\C}_{ij}\subseteq Q_{ij}$ for all $i,j\in\{1,2\}$, therefore
the theorem follows from the following stronger lemma.
\end{proof}

\begin{lemma}
If $x_1\in Q_{i_1i_2}, x_2\in Q_{i_2i_3}, \dots, x_n\in
Q_{i_ni_{n+1}}$ (where $i_1,\dots,i_{n+1}\in\{1,2\})$ then
$$\phi(x_1 x_2 \cdots x_n)=0.$$
If furthermore $(i_1,i_{n+1})\neq (1,2)$ then
$$E(x_1 x_2 \cdots x_n)=0.$$
\end{lemma}

\begin{proof}
We shall prove the lemma by induction over $n$.

Since $Q_{12}=\F_0\oplus \{x\in Q_{12}:E(x)=0 \}$, therefore we can assume
that for all $k$ such that $(i_k,i_{k+1})=(1,2)$ we have either $x_k\in\F_0$
or $E(x_k)=0$.

Note that if $(i,j)\neq (1,2)$ then $Q_{ij}\subseteq P^{\F}_{ij}$.
(For $(i,j)=(2,1)$ observe that if $x=b_1
d_1 b_2 d_2 \cdots d_{n-1} b_n$, where
$b_1,\dots,b_n\in\B_0$ and $d_1,\dots,d_{n-1}\in\D_0$ then for every
$d\in\D$ we have
$$\phi[E(x) d]=\phi[E(x d)]=
\phi(x d)=\phi\big[x \big(d-\phi(d)\big)\big]+\phi(x) \phi(d)=0$$ by
freeness. Since $\phi$ is faithful we have that $E(x)=0$.)

Furthermore $\{x\in Q_{12}:E(x)=0\}\subseteq P^{\F}_{12}$;
therefore if $x_1\in Q_{i_1i_2}, x_2\in
Q_{i_2i_3}, \dots, x_n\in Q_{i_ni_{n+1}}$ have the additional property that
$E(x_k)=0$ for all $k$ such that $(i_k,i_{k+1})=(1,2)$ then the condition
\ref{war:momenty} of Theorem \ref{thedef:centralne} implies
$E(x_1 x_2\dots x_n)=0.$

If the above case does not hold, let $k$ the index with the property
that $(i_k,i_{k+1})=(1,2)$ and $x_k\in\F_0$.
If $k=1$ then $E(x_1 x_2 \dots x_n)=x_1 E(x_2 \dots x_n)=0$ by inductive
hypothesis, similar argument is valid for $k=n$.

For $2\leq k\leq n-1$ let us consider a product
$$x_1 x_2 \cdots x_{k-2} (x_{k-1} x_k x_{k+1}) x_{k+2} \cdots x_n.$$
We will show that $x_{k-1} x_k  x_{k+1}\in
Q_{i_{k-1}i_{k+2}}$, therefore the inductive hypothesis can be applied.

For $(i_{k-1},i_{k+2})\neq (1,2)$ the proof is very simple.
For  $(i_{k-1},i_{k+2})=(1,2)$ we use the assumption that $\phi$ is tracial
and that $\B$ and $\F$ are free.
\end{proof}

\subsection {Matrix algebra}
\label{sect:macierze}
We will consider here an important example of a pair of compatible
structures on the matrix algebra.

Let us fix $N\in\N$ and
let $(\tilde{\A},\tilde{\phi})$ be a scalar probability space.
We define $\A=\M_N(\C)\otimes\tilde{\A}$ to be the matrix algebra over
$\tilde{\A}$.
We have that $(\M_N(\C)\subseteq\A,E)$ is an $\M_N(\C)$--valued probability space,
where the expectation value
$E:\A\rightarrow\M_N(\C)$ is defined by $E=\rm{Id}\otimes \tilde{\phi}$, i.e.
$$E\big( [ a_{ij}]_{1\leq i,j\leq N} \big)=[ \tilde{\phi}(a_{ij}) ]_{1\leq i,j \leq N}.$$

We can furnish $\A$ with a structure of a scalar probability space $(\A,\phi)$,
where $\phi(  [ a_{ij}]_{1\leq i,j\leq N})=\frac{1}{N} \sum_k \tilde{\phi}(a_{kk})$.

We will denote by $\tilde{\kk}^{\C}$ (resp.\ $\kk^{\C}$) the scalar cumulant
function in the algebra $(\tilde{\A},\tilde{\phi})$ (resp.\ $(\A,\phi)$).

The following theorem goes back to the work of Nica, Shlyakhtenko, and Speicher
and can be found in \cite{SpIHP}.
\begin{theorem}  \label{theo:kiedywolneodskalarow}
An element $x=[x_{ij}]_{1\leq i,j\leq N}\in\A$ and the algebra $\M_N(\C)\subseteq
\A$ are free (with respect to the state $\phi$)
 iff for every $n\geq 1$ there
exists $c_n\in \C$ such that
$$ \tilde{\kk}^{\C}[ x_{i_1i_2}\otimes x_{i_2i_3}\otimes\cdots\otimes
x_{i_ni_1}]=c_n$$
for every $i_1,\dots,i_n\in\{1,\dots,n\}$
and all cumulants
$\tilde{\kk}^{\C}[ x_{i_1j_1}\otimes x_{i_2j_2}\otimes\cdots\otimes
x_{i_nj_n}]$ that are not of this type vanish.

If the above holds we have that
$$ c_n=N^{n-1} \kk^{\C}(\underbrace{x\otimes\cdots\otimes x}_{n\ {\rm
times}}). $$
\end{theorem}

The following result is due to Nica, Shlyakhtenko, and Speicher and will be
published in \cite{NiShSp2}.
In order to be self-contained we indicate the proof.
\begin{lemma} \label{lem:kumulantyoperatorowe}
The cumulants $\kk^{\M_N(\C)}$ in the probability space
$(\M_N(\C)\subseteq \A,E)$ and the cumulants $\tilde{\kk}^{\C}$
in the probability space $(\tilde{A},\tilde{\phi})$ are related as follows:
\begin{multline}
\kk^{\M_N(\C)}[(m_1\otimes z_1)\otimes
\cdots\otimes (m_n\otimes z_n)]=  \\
 \tilde{\kk}^{\C}(z_1\otimes\cdots\otimes z_n) (m_1 m_2 \cdots m_n)
\label{eq:wzor}
\end{multline}
for $m_1\otimes z_1,\dots,m_n\otimes z_n\in \M_N(\C)\otimes \tilde{A}$.
\end{lemma}
\begin{proof}
We can take (\ref{eq:wzor}) as a definition of the function
$\kk^{\M_N(\C)}$. It is not difficult to show that the
corresponding compound functions are given by
$$\kk_{\pi}^{\M_N(\C)}[(m_1\otimes z_1)\otimes\cdots\otimes
(m_n\otimes z_n)]=\tilde{\kk}_{\pi}^{\C}(z_1\otimes\cdots\otimes z_n)
(m_1 m_2 \cdots m_n) $$
for every $\pi\in\NC(\{1,\dots,n\})$. Thus we see that the defining relation
(\ref{eq:momentyprzezkumulanty}) for cumulants is indeed
fulfilled.
\end{proof}

\subsection [Application: upper triangular matrices\dots]{Application:
upper triangular matrices of Dykema and Haagerup}
In this section we present an application of the tools of the
operator--valued free probability, namely using significantly simpler methods
we prove that circular free Poisson elements have a certain upper triangular
representation.

We say that an element $x$ of a scalar probability space $(\tilde{\A},
\tilde{\phi})$ is
a circular free Poisson element \index{circular free Poisson element}
with parameter $c\geq 1$ if $x$ is $R$--diagonal and if for every $n\in\N$,
$n\geq 1$ we have
$$\kk (\underbrace{x x\gwia \otimes x x\gwia\otimes \cdots\otimes
 x x\gwia}_{n\ {\rm times}})=\kk(\underbrace{x\gwia x \otimes x\gwia
x\otimes\cdots\otimes x\gwia x}_{n\ {\rm times}})=c. $$

We say that an element $x$ of a scalar probability space
$(\tilde{\A},\tilde{\phi})$ is
circular \index{circular element}
if $x$ is a free Poisson element with parameter $c=1$.

The following lemma follows directly from Theorem \ref{theo:krawczykspeicher}
 of Krawczyk and Speicher.
\begin{lemma} \label{lem:lematspeicherowy}
If elements $x_1,x_2,\dots,x_{2n}\in(\A,\phi)$
have the property that for every odd number
$2k+1$ and all indices $i_1,\dots,i_{2k+1}\in\{1,\dots,2n\}$ we have that
$$\kk(x_{i_1}\otimes\cdots\otimes x_{i_{2k+1}})=0$$
then
$$\kk (x_2 x_3\otimes x_4 x_5\otimes\cdots\otimes x_{2n-2}
x_{2n-1}\otimes x_{2n} x_{1})=\sum_{\pi} \kk_{\pi}(x_1 \otimes x_2\otimes
\cdots \otimes x_{2n}),$$ where the sum is taken over all noncrossing
partitions $\pi$ of the set $\{1,2,\dots,2n\}$ such that $\big\{
\{1,2\},\{3,4\},\dots,\{2n-1,2n\} \big\} \prec \pi$. \end{lemma}

\begin{theorem} Let $c\geq 1$ and $N\in\N$, and let $(\tilde{\A},
\tilde{\phi})$ be a
scalar probability space with random variables
$a_1,\dots,a_N\in\tilde\A$ and $b_{ij}\in\tilde\A$ ($1\leq i<j\leq N$), where $a_j$ is a
circular free Poisson element with parameter $(c-1)N+j$ and each
$b_{ij}$ is a circular element,
and where the family
$$\big( (\{a_j\gwia,a_j\})_{1\leq j\leq N}, (\{b_{ij}\gwia,b_{ij}\})_{1\leq
i<j \leq N} \big) $$ is free. Consider the scalar
probability space $(\A,\phi)$ defined as in Sect.\ \ref{sect:macierze}
of this Section  and consider the random variable
\begin{equation} x=\frac{1}{\sqrt{N}} \left[ \begin{array}{cccccc}
a_1    & b_{12} & b_{13} & \cdots & b_{1,N-1} & b_{1,N} \\
0      & a_2    & b_{23} & \cdots & b_{2,N-1} & b_{2,N} \\
0      & 0      & a_3    & \ddots & \vdots    & b_{3,N} \\
\vdots & \vdots & \ddots & \ddots & \ddots    & \vdots  \\
0      & 0      & \cdots & 0      & a_{N-1}   & b_{N-1,N} \\
0      & 0      & \cdots & 0      & 0         & a_N
\end{array} \right]. \label{eq:haadyk} \end{equation}
Then $x\in(\A,\phi)$ is a circular free Poisson element of parameter $c$.
\end{theorem}
\begin{proof}
Denote the entries of the matrix $\sqrt{N} x$ by $y_{ij}$ ($1\leq i,j\leq N$).
We use Lemma \ref{lem:lematspeicherowy} to compute
$$\kk[y_{i_1 j_1} (y_{k_1 j_1})\gwia \otimes y_{i_2 j_2} (y_{k_2 j_2})\gwia
\otimes \cdots \otimes y_{i_n j_n} (y_{k_n j_n})\gwia].$$
Since the family $\{y_{ij}\}$ is free, the only nonvanishing cumulants are
of the type
\begin{multline} \kk[y_{i_1 a} (y_{i_2 a})\gwia \otimes y_{i_2 a} (y_{i_3 a})\gwia
\otimes \cdots \otimes y_{i_n a} (y_{i_1 a})\gwia]=\\
 \sum_{\pi} \kk_{\pi}[(y_{i_1 a})\gwia
\otimes y_{i_1 a}\otimes \cdots \otimes (y_{i_n a})\otimes \gwia y_{i_n
a}], \label{eq:rownanie1500}
\end{multline}
where the sum is taken over such $\pi$ as in Lemma
\ref{lem:lematspeicherowy}. This expression will be evaluated in
Lemma \ref{lem:lematobliczeniowy} below. It follows that the only
nonvanishing cumulants of the matrix $z=xx\gwia$ are $$\kk(z_{i_1
i_2}\otimes \cdots \otimes z_{i_{n-1} i_n} \otimes z_{i_n
i_1})=\frac{c}{N^{n-1}}.$$

Let the expectation value $E:\A\rightarrow\M_N(\C)$
be as in Sect.\ \ref{sect:macierze}.
Condition 2 of Theorem
\ref{thedef:centralne} and Lemma
\ref{lem:kumulantyoperatorowe}
show that $x$
is $R$--diagonal with amalgamation over $\M_N(\C)$.
Furthermore, the above computation and Theorem
\ref{theo:kiedywolneodskalarow} show that $xx\gwia$ and $\M_N(\C)$ are free.
It follows from Theorem \ref{theo:rdiagonalnosc} that $x\in(\A,\phi)$ is
$R$--diagonal.

The distribution of any $R$--diagonal element $t$ is uniquely determined by the
distribution of $t t\gwia$ and
the above calculation shows that the distribution of $xx\gwia$ coincides with
a distribution of $s s\gwia$, where $s$ is a free circular
Poisson element with parameter $c$.
\end{proof}
\begin{remark}
We have shown that the matrix (\ref{eq:haadyk}) is $R$--diagonal with
amalgamation over $\M_N(\C)$ and that $xx\gwia$ and $\M_N(\C)$ are free;
it coincides therefore with the upper triangular form we will construct
in Theorem \ref{theo:utmform}.
\end{remark}
\begin{lemma} \label{lem:lematobliczeniowy} Let $(y_{ij})$ be as above.
For every $1\leq a,i_1,i_2,\dots,i_n\leq N$ we have
\begin{multline} \sum_{\pi} \kk_{\pi} [(y_{i_1 a})\gwia \otimes y_{i_1 a}
\otimes (y_{i_2 a})\gwia \otimes y_{i_2 a} \otimes \cdots \otimes
(y_{i_n a})\gwia \otimes y_{i_n a}]\label{eq:rownanie998} \\
=\left\{ \begin{array}{cl} (c-1)N+a  & {\rm if}\
\max(i_1,i_2,\dots,i_n)=a \\
1 & {\rm if}\ \max(i_1,i_2,\dots,i_n)<a \\
0 & {\rm if}\ \max(i_1,i_2,\dots,i_n)>a
\end{array} \right. , \end{multline}
where the sum is taken over all $\pi\in{\rm NC}(\{1,\dots,2n\})$ such that
$\big\{\{1,2\},\{3,4\},\dots,\{2n-1,2n\} \big\}\prec\pi$.
\end{lemma}
\begin{proof}
If $\max(i_1,\dots,i_n)>a$ then the statement of the lemma follows
trivially since one of the factors is equal to zero.

If $i_1=\cdots=i_n$ then the statement follows easily from Lemma
\ref{lem:lematspeicherowy}.

For the other cases we will use the induction with respect to
$\max(i_1,\dots,i_n) - \min(i_1,\dots,i_n)$.
We define $s=\min(i_1,\dots,i_n)$,
$B=\{2m-1,2m:i_m=s \}\subseteq\{1,2,\dots,2n\}$ and 
$A=\{1,2,\dots,2n\}\setminus B=\{2m-1,2m:i_m>s\}$.

For every noncrossing partition
$$\sigma=\bigl\{ \{i_{1,1},i_{1,2},\dots,i_{1,m_1} \},
\{i_{2,1},i_{2,2},\dots,i_{2,m_2} \}, \dots, 
\{i_{k,1},i_{k,2},\dots,i_{k,m_k} \} \bigr\}$$
 of any subset of the set $\{1,2,\dots,2n\}$ (where 
$i_{j,1}<i_{j,2}<\cdots<i_{j,m_j}$ for every value of $j$)
we define 
$$\kk_{\pi}( a_1,\dots,a_{2n})=
\prod_{1\leq l\leq k} \kk(a_{i_{l,1}}\otimes
a_{i_{l,2}}\otimes\cdots\otimes a_{i_{l,m_l}}).$$

Due to the freeness every partition $\pi$ which contributes
nontrivially in the sum
(\ref{eq:rownanie998}) can be written as $\pi=\pi_A\cup \pi_B$, where
$\pi_A\in\NC(A)$ and $\pi_B\in\NC(B)$; therefore the left--hand side of
(\ref{eq:rownanie998}) is equal to 
\begin{multline*} \sum_{\pi_A} \Big[
\kk_{\pi_A} [(y_{i_1 a})\gwia,
y_{i_1 a},\dots,(y_{i_n a})\gwia, y_{i_n a}] \times \\
 \sum_{\pi_B} \kk_{\pi_B} [(y_{i_1 a})\gwia, 
y_{i_1 a},\dots,(y_{i_n a})\gwia ,y_{i_n a}] \Big], \end{multline*}
where the first sum is taken over $\pi_A\in\NC(A)$ such that
$\big\{ \{2m-1,2m\} :i_m>s \big\}\prec \pi$
and the second sum is taken over $\pi_B\in\NC(B)$ such that
$\big\{ \{2m-1,2m\} :i_m=s \big\}\prec
\pi_B\prec(\pi_A)^{c}$.

For every $\pi_A$ we have that $\sum_{\pi_B} \kk_{\pi_B} [(y_{i_1 a})\gwia, 
y_{i_1 a},\dots,(y_{i_n a})\gwia ,y_{i_n a}]$ is equal to a product of 
$|(\pi_A)^{c}|$ expressions of type $\kk[y_{s a} y_{s
a}\gwia \otimes y_{s a} y_{s a}\gwia \otimes \cdots \otimes y_{s a} y_{s
a}\gwia ] $ which is always equal to $1$ (we use here that
$s=\min(i_1,\dots,i_n)<\max(i_1,\dots,i_n)\leq a$, so $s<a$ and $y_{s a}$ is
an element from above of the diagonal of the matrix).

Therefore we showed that the left hand side of (\ref{eq:rownanie998}) is
equal to $$\sum_{\pi_A} \kk_{\pi_A} [(y_{i_1 a})\gwia,
y_{i_1 a},\dots,(y_{i_n a})\gwia, y_{i_n a}]$$ and the inductive hypothesis
can be applied.
\end{proof}


\subsection {Generalized matrix algebra}
\label{sect:matrix}

We will consider a generalization of the example from Sect.\
\ref{sect:macierze}
in which the matrix structure of $\A$ comes from a sequence of projections.

At the beginning we will show how to recover scalar structure
from the compatible operator--valued one and vice versa.

\subsubsection {Scalar structure induces the operator--valued structure}
\label{subsect:scalarinduces}
Let a scalar probability space $(\A,\phi)$ be given and
suppose that there are orthogonal projections $e_1,\dots,e_n\in \A$
such that $e_1+\cdots+e_n=I$.

Let $\D$ be the unital algebra generated by $e_1,\dots,e_n$. We define
an expectation value $E:\A\rightarrow\D$ by
$$E(z)=\sum_{k} \frac{\phi(z e_k)}{\phi(e_k)} e_k.$$
In this way we have furnished
$\A$ with a compatible structure of a $\D$--valued operator space.

\subsubsection {Operator--valued structure induces the scalar structure}
\label{subsect:operatorinduces}
Let an operator--valued probability space $(\D\subseteq\A,E)$ be given,
such that $\D$ is a unital algebra generated by orthogonal projections
$e_1,\dots,e_n$ and let $t_1,\dots,t_n>0$, $t_1+\cdots+t_n=1$.

We define $\phi:\A\rightarrow\C$ by $\phi(z)=\phi_0[E(z)]$, where
$\phi_0:\D\rightarrow\C$ is a linear functional defined by
$\phi_0(e_i)=t_i$. In this way we have furnished $\A$ with the structure of a
scalar probability space $(\A,\phi)$ and
we have that $\phi(z)=\phi[E(z)]$ for every $z\in\A$.

\subsubsection {Scalar random variables}
It is not difficult to see that in this example every scalar random
variable $z$ is a diagonal matrix of the form
$$z=\left[ \begin{array}{ccccc}
e_1 z e_1 & 0 & \cdots &  & 0 \\
0 & e_2 z e_2 &        &  &   \\
\vdots & & \ddots     &  & \vdots \\
  &     &  & e_{n-1} z e_{n-1} & 0 \\
0 & &\cdots  &  0 & e_n z e_n \end{array} \right]. $$

\subsubsection {Upper triangular square root}
\begin{lemma} \label{lem:pierwiastek}
Let $\A\subseteq B(\Ha)$ be a von Neumann algebra and
let $e_1,\dots,e_n\in\A$ be orthogonal projections such that
$e_1+\cdots+e_n=I$. If $y\in\A$ is a strictly positive
operator then there
exists an operator $x\in\A$ such that $y=x x\gwia$ and furthermore
$x$ is an upper triangular matrix (\ref{eq:utm1}) with respect to
$e_1,\dots,e_n$. 
\end{lemma} 
\begin{proof} 
We consider first the case $n=2$. 
For $$x=\left[ \begin{array}{cc} x_{11} & x_{12} \\ 0 & x_{22} \end{array}
\right],$$ where $x_{ij}=e_i x e_j$,
the equation $x x\gwia=y$ reads
$$ \left[ \begin{array}{cc} x_{11} x_{11}\gwia+x_{12} x_{12}\gwia & 
x_{12} x_{22}\gwia \\ x_{22} x_{12}\gwia & x_{22} x_{22}\gwia
\end{array} \right]=\left[ \begin{array}{cc} y_{11} & y_{12} \\ y_{21} &
y_{22} \end{array} \right],$$
where $y_{ij}=e_i y e_j$. 
We solve this system of equations and find $x_{22}=(y_{22})^{\frac{1}{2}}$,
$x_{12}=y_{12} x_{22}^{-1}$ (since $y$ is strictly positive,
$x_{22}$ is invertible).

The equation
$x_{11} x_{11}\gwia=y_{11}-x_{12} x_{12}\gwia$
has a solution since 
the operator on the right hand side is strictly positive; for every $v\in e_1
\Ha$, $v\neq 0$  we have
\begin{multline*} \langle v, \bigl[ y_{11} -y_{12}
(y_{22}^{-1}) y_{12}\gwia\bigr] v \rangle= \\ \left\langle
\left[ \begin{array}{c} v \\ -y_{22}^{-1} y_{12}\gwia v \end{array}
\right],\left[ \begin{array}{cc} y_{11} & y_{12} \\ y_{21} & y_{22}
\end{array} \right] \left[ \begin{array}{c} v \\ -y_{22}^{-1}
y_{12}\gwia v \end{array} \right] \right\rangle> 0.
\end{multline*}

The general case follows from this special case $n=2$ by induction as
follows. We apply the lemma to the operator $y$ and two
projections $f_1=e_1+\cdots+e_{n-1}$ and $f_2=e_n$ and thus get
an operator $x_0$ such that $y=x_0 x_0\gwia$ and
$$x_0=\left[ \begin{array}{cc} f_1 x_0 f_1 & f_1 x_0 f_2 \\ 0 & 
f_2 x_0 f_2 \end{array} \right].$$
By induction hypothesis, we can now apply the lemma for the operator
$$y_1=(f_1 x_0 f_1) (f_1 x_0 f_1)\gwia\in B[f_1 \Ha]$$
and the projections $e_1,\dots, e_{n-1}$ and thus get an $x_1\in B[(1-e_n)\Ha]$
such that $x_1 x_1\gwia=y_1$ and $x_1$ is an upper triangular matrix with 
respect to $e_1,\dots,e_{n-1}$. 
The element $$x=\left[ \begin{array}{cc} x_1 & f_1 x_0 f_2 \\ 
0 & f_2 x_0 f_2 \end{array} \right] $$ 
gives the wanted triangular square root with respect to $e_1,\dots,e_n$. 
\end{proof}


\section {The main result}
\label{chap:main}
\subsection {Upper triangular form for $R$--diagonal elements}
\begin{theorem}
\label{theo:utmform}
Every $R$--diagonal
operator has the upper triangular form property.
\end{theorem} \begin{proof}
Let $x'\in(\A',\phi')$ be an $R$--diagonal operator.
Our goal is to construct a scalar probability space $(\A,\phi)$ and
elements $x,e_1,\dots,e_N\in\A$ such that
$e_1,\dots,e_N$ are orthogonal projections in respect to which
$x$ is in the upper triangular form and
such that $x\in(\A,\phi)$ has the same distribution as $x'$. 
If the condition on spectra (\ref{eq:spektra2}) holds then 
it follows that also the original element $x'$ has the upper triangular
property (cf \cite{DyHa2000}). 

First of all note that, for every $s\geq 0$, the subspace
$V_s$ defined in (\ref{eq:niezmienniczaHa}) fulfills
$$V_s=\overline{ \left\{v\in \Ha: \limsup_{n\rightarrow\infty}
\sqrt[2n]{\langle v,(x\gwia)^n x^n v\rangle}\leq s \right\} }, $$
therefore $p_s\in\left\{ x\gwia x, (x\gwia)^2 x^2, \dots \right\}''$, hence
$p_s$ and $x x\gwia$ are free (Theorem \ref{theo:nowardiagonalnosc}).
The projections $e_i$ that we want to construct will turn out at the end
of the proof to be equal to $e_i=p_{s_i}-p_{s_{i-1}}$; these are heuristical
motivations why we will choose $e_1,\dots,e_n$ to be free from $x x\gwia$.

{\bf Construction of $x$.}
Let us first restrict to the case where $x'(x')\gwia$ is strictly positive.

For given $t_1,t_2,\dots,t_N>0$, $t_1+\cdots+t_N=1$ let us consider an
enlargement $(\A_1,\phi_1)$ of the probability space
$\Big({\rm Alg}\big(I,x' (x')\gwia\big),\phi'\Big)$.
This enlargement should contain
orthogonal projections $e_1,e_2,\dots,e_N\in\A_1$ such that $e_1+\cdots+e_N=I$,
$\phi_1(e_i)=t_i$ and such that $x' (x')\gwia$ and $\{e_1,\dots,e_N\}$ are
free (this enlargement is given by the Avitzour--Voiculescu
free product of algebras \cite{Avi,Voi1985}).
\index{free product of algebras}

Lemma \ref{lem:pierwiastek} gives us  an element $z\in\A_1$
such that $z z\gwia=x' (x')\gwia$ and such that $z$ has the upper
triangular form (\ref{eq:utm1}).

We define the operator--valued structure $(\D\subseteq\A_1,E_1)$ on $\A_1$ as
in Sect.\ \ref{subsect:scalarinduces} of Section  \ref{chap:compatible}.
The free product of algebras with amalgamation over $\D$ \cite{VDN}
\index{free product of algebras with amalgamation} gives
us an enlargement $(\D\subseteq\A,E)$ of the space $(\D\subseteq\A_1,E_1)$ and
an operator $u\in\A$ which is a scalar Haar unitary such that $\{u,u\gwia\}$ and
$\A_1$ are free
with amalgamation over $\D$. As in Sect.\ \ref{subsect:operatorinduces}
of Section  \ref{chap:compatible}
we recover the scalar--valued structure $(\A,\phi)$ on $\A$.
Note that $(\A,\phi)$ is an enlargement of $(\A_1,\phi_1)$.

Theorem \ref{thedef:centralne} gives us that $zu$ is
$R$--diagonal with amalgamation over $\D$; we know that $(zu)(zu)\gwia=x'
(x')\gwia$ and $\D$ are free, therefore Theorem \ref{theo:rdiagonalnosc}
assures us that $zu$ is $R$--diagonal. The distribution of any
$R$--diagonal element $y$ is uniquely determined by the distribution of $y
y\gwia$; therefore the elements $x'$ and $zu$ are identically distributed.

Since $u$ is a diagonal matrix, the operator $x:=zu$ is the wanted upper
triangular matrix.

{\bf Estimates for spectra.}

The set $\sigma_i$, defined to be the spectrum of $e_i x e_i\in e_i \A e_i$ is 
contained in the spectrum of
$p x p\in p \A p$, where $p=e_i+e_{i+1}+\cdots+e_N$. 
 
We can furnish $p \A p$ with the operator--valued structure $(\D'\subseteq 
p \A p,E')$ in the standard way (with respect to orthogonal projections
$e_i,\dots,e_N$); note that $\D'=p \D p$ and $E'(y)=E(y)$
for every $y\in p\A p$. For $z_1,\dots,z_k\in p\A p$ we have
$$\kk^{\D}(z_1\otimes\cdots\otimes z_k)=
\kk^{\D}(p z_1\otimes\cdots\otimes z_k p)=
p \kk^{\D}(z_1\otimes\cdots\otimes z_k) p $$
and therefore $\kk^{\D}(z_1\otimes\cdots\otimes z_k)\in \D'$.
This observation implies that the (restriction of) 
cumulant function $\kk^{\D}$ in
the probability space $(\D\subseteq \A,E)$ fulfills the 
definition of the cumulant function $\kk^{\D'}$ in the space $(\D'\subseteq
\A',E')$ so
$$\kk^{\D'}(z_1\otimes \cdots \otimes z_k)=\kk^{\D}(z_1\otimes \cdots 
\otimes z_k).$$

We apply now Theorem \ref{theo:wolneodskalarow2} and see that 
$(p x p)(p x p)\gwia=
p x x\gwia p$ and $\D'$ are free in the scalar probability space
$(p A p,\frac{1}{\phi(p)} \phi)$. We also see that 
$p x p\in(\D'\subset p\A_2 p,E_2')$
is $R$--diagonal with amalgamation over $\D'$. Theorem \ref{theo:rdiagonalnosc}
assures us that $p x p\in(p
A p,\frac{1}{\phi(p)} \phi)$ is $R$--diagonal.

Lemma \ref{lem:stransformataipromienspektralny} implies that 
$$\sigma_i\subseteq\left\{z\in\C: |z|\geq R_i\right\},$$ 
where $$R_i= \frac{1}{S^{p \A p}_{(p x p) (p x p)\gwia}(-1)}$$ 
 
Due to the upper triangular form of $x$ we have $(p x p) (p x 
p)\gwia=p xx{}\gwia p$ and Lemma \ref{lem:stransformatakompresji}
implies that 
\begin{equation} \label{eq:wzornaR} 
R_i=\frac{1}{S^{\A}_{x x\gwia}[-(t_i+t_{i+1}+\cdots+t_N)]} 
\end{equation}

Similarly as above, we see that the spectrum of
$e_{i-1} x e_{i-1}\in e_{i-1} \A e_{i-1}$ is a subset of
the spectrum of $q x q\in q \A q$, where $q=I-p=e_1+\cdots+e_{i-1}$.

Lemma \ref{lem:stransformataipromienspektralny} implies that
$$\sigma_{i-1}\subseteq\left\{z\in\C: |z|\leq r_i\right\},$$
where $$r_i=S^{q\A q}_{[(qxq)(qxq)\gwia]^{-1}}(-1)=
\phi[(q x q) (q x q)\gwia],$$
where $[(qxq)(qxq)\gwia]^{-1}$ denotes the inverse of $(qxq)(qxq)\gwia\in
qAq$.
 
Let us introduce the notation 
$$x =\left[ \begin{array}{cc} q x q & q x p \\ 
 0 & p x p 
\end{array} \right]=\left[\begin{array}{cc} a & b \\ 
0 &c \end{array}\right],$$ 
therefore
$$x^{-1}=\left[ \begin{array}{cc} a^{-1} & -a^{-1} b c^{-1} \\
0 & c^{-1} \end{array} \right] $$ and 
$$(x x\gwia)^{-1}= 
\left[ \begin{array}{cc} (a^{-1})\gwia a^{-1} & -(a^{-1})\gwia a^{-1} b c^{-1} \\
-(c^{-1})\gwia b\gwia (a^{-1})\gwia a^{-1} & (c^{-1})\gwia b\gwia (a^{-1})\gwia
a^{-1} b c^{-1}+(c^{-1})\gwia c^{-1} \end{array} \right]$$
and therefore
\begin{equation} [(q x q)^{-1}]\gwia (q x q)^{-1}=
q (x x\gwia)^{-1} q, \label{eq:dziwnywzornaodwrotnosc}
\end{equation} 
is well--defined, where $(q x q)^{-1}$ denotes the inverse of $q x
q\in q \A q$. As before show that $(q x q)^{-1}\in(q A
q,\frac{1}{\phi (q)} \phi)$ is $R$--diagonal and then 
Eq.\ (\ref{eq:dziwnywzornaodwrotnosc}) and Lemma
\ref{lem:stransformatakompresji} imply that
\begin{equation} \label{eq:wzornar} 
r_i= 
S^{\A_2}_{(x x\gwia)^{-1}}[-(t_1+\cdots+t_{i-1})]. 
\end{equation} 
 
Lastly, Lemma \ref{lem:stransformataodwrotnego} shows that $R_i=r_i$.
 
Proposition 2.2 from \cite{DyHa2000} shows that
$e_1+e_2+\cdots+e_i=p_{r_i}$, where $p_{s}$ is the projection 
on the space (\ref{eq:niezmienniczaHa}).

{\bf The general case.}
In the general case where $x'(x')\gwia$ might have a kernel,
for every $\epsilon>0$ we construct as above an $R$--diagonal element
$x_{\epsilon}\in(\A,\phi)$ in the triangular form with respect to
projections $e_1,\dots,e_N\in\A$, such that $x_{\epsilon} x_{\epsilon}\gwia$
and $x' (x')\gwia+\epsilon$ are identically distributed.
 
There exists a sequence $(\epsilon_n)$ of positive numbers which converges
to zero such that the sequence
$\phi_1[ P(x_{\epsilon}, x_{\epsilon}\gwia, e_1, e_2, \dots, e_N)]$ converges 
for every polynomial $P$ of $N+2$ noncommuting variables. 
We define $\tilde{\A}$ to be the free algebra generated by elements 
$x,x\gwia,e_1,\dots,e_N$ and a state $\tilde{\phi}$ on $\tilde{\A}$ by 
$$\tilde{\phi}[ P(x,x\gwia,e_1,\dots,e_N) ]= \lim_{n \rightarrow \infty} 
\phi[ P(x_{\epsilon_n},x_{\epsilon_n}\gwia,e_1,\dots,e_N) ]. $$
 
We have that $x$ is in the upper triangular form and $x$ and $x'$ are
identically distributed. 

The spectra estimates can be done similarly.
\end{proof}

\begin{remark}
\label{remark:equality}
In the case $N=2$ our proof shows that we can write a given
$R$--diagonal distribution in the form of a $2\times 2$ matrix
$$x=\begin{pmatrix} qxq&qxp\\ 0&pxp \end{pmatrix}=
\begin{pmatrix} a&b\\ 0&c \end{pmatrix},$$
where the spectrum and the distribution of the diagonal elements 
$a$ and $c$ can be described quite concretely.
 
We calculated the spectra as 
$$\sigma(a)=
\{z\in\C\mid R_0\leq\vert z\vert\leq R_{1}\},
\qquad \sigma(c)=\{z\in\C\mid R_1\leq\vert z\vert\leq R_2\},$$
 where
\begin{align*}
R_0&=\frac 1{S_{xx\gwia }(-1)}\\
R_1&=S_{(xx\gwia)^{-1} }(-t)=\frac 1{S_{xx\gwia }(t-1)}\\
R_2&=S_{(xx\gwia)^{-1} }(-1).
\end{align*}
with $t=\phi(q)$. 

The distribution of $c$ is determined by the facts that
$c$ is $R$--diagonal and that $cc\gwia $ has the 
same distribution as $pxx\gwia p$. The distribution of $a$ is determined 
by the facts that $a$ is $R$--diagonal and that $(a^{-1})\gwia a^{-1}$ has the same distribution as 
$q(xx\gwia )^{-1}q$. These requirements determine both distributions 
uniquely since we also have by construction that $p$ is free from $xx\gwia $. 
(In the latter case one should also note that, by
\cite{HL}, $R$--diagonality of $a$ implies $R$-diagonality of
$a^{-1}$.) 
 
Since the upper triangular form realization is unique, one can get
the realizations of $x$ as triangular $N\times N$-matrices by
iteration of the case for $N=2$. But this implies that also in the
case of general $N$ the distribution of all diagonal elements can
be calculated concretely by iteration of the above two
prescriptions and that the spectra in this case are given by
$$\sigma(e_kxe_k)=\{z\in\C\mid R_{k-1}\leq\vert z\vert\leq
R_{k}\},$$ where $$R_k=S_{(xx\gwia)^{-1} }(-t)=\frac
1{S_{xx\gwia }(t-1)}\qquad\text{where $t:=\phi(e_1+\dots+e_{k-1})$.}$$ 
\end{remark}

\subsection {Brown measure of $R$--diagonal elements}
\label{sec:brownmeasure}
\index{Brown measure of $R$--diagonal elements}

If $x\in(\A,\phi)$ is $R$--diagonal then for every complex number
$c$ such that $|c|=1$ we have that elements $x$ and $cx$ are
identically distributed (Theorem \ref{thedef:centralne}, condition
\ref{war:kazdyunitarny}), therefore the Brown measure of $x$ is
rotationally invariant. Eq.\ (\ref{eq:wzornaR}),
(\ref{eq:wzornar}) and the discussion from Sect.\
\ref{subsect:miarabrowna} of Section  \ref{chap:introduction}
imply the following result, which was obtained for the first time by
Larsen and Haagerup \cite{HL,La}.

\begin{theorem}
Let $x\in (\A,\phi)$ be $R$--diagonal. Brown spectral
distribution measure
$\mu_x$ of $x$ is the rotationally invariant measure
on $\C=[0,\infty)\times \T$
(these sets are identified by the polar coordinates
$(r,\psi)\mapsto r e^{i \psi}$)
$$\mu_x=\nu_x \times \chi,$$
where $\nu_x$ is the probability measure on $[0,\infty)$ with the
inverse of the distribution function $F_x(s)=\nu_x[0,s]$ given by 
$$F_x^{\langle -1\rangle}(t)= 
\frac{1}{S_{x x\gwia}(t-1)}=S_{(x x\gwia)^{-1}}(-t) 
\qquad\mbox{for } 0\leq t\leq 1.$$
\end{theorem}
 
\section [Upper triangular form for products]{Upper
triangular form for products of free $R$--diagonal elements}
\label{chap:produkt}
Let $x_1',x_2'\in(\A',\phi')$
be $R$--diagonal, let $\{x_1',(x_1')\gwia \}$ and $\{x_2',(x_2')\gwia \}$
be free and let $t_1,\dots,t_N>0$, $t_1+\cdots+t_N=1$ be given. 
 
For $i\in\{1,2\}$
let $x_i,e_1,e_2,\dots,e_n\in (\A_i,\phi_i)$ be operators given by
Theorem \ref{theo:utmform} applied for $x_i'$; we furnish $\A_i$ with
a structure of an operator--valued space $(\D\subseteq\A_i,E_i)$ as
in Sect.\ \ref{subsect:scalarinduces} of Section  \ref{chap:compatible}.

We take now the free product
of the algebras $(\D\subseteq\A_i,E_i)$, $i\in\{1,2\}$ with amalgamation over
$\D$ and denote the resulting space by $(\D\subseteq\A,E)$.
As in Sect.\ \ref{subsect:operatorinduces} of Section  \ref{chap:compatible}
we recover the scalar--valued structure $(\A,\phi)$ on $\A$.

\begin{theorem}
The products $x_1 x_2\in (\A,\phi)$ and $x_1' x_2'\in (\A',\phi')$ are
identically distributed.
\end{theorem}
\begin{proof}
It is well-known \cite{NiSp1997A} that $x_1' x_2'\in (\A',\phi')$ is $R$--diagonal (this follows e.g.,
by condition \ref{war:kazdyunitarny} of Theorem
\ref{thedef:centralne}).
We will show that $x_1 x_2\in (\A,\phi)$ is also $R$--diagonal.
Firstly, note that $x_1 x_2\in (\D\subseteq \A,E)$ is $R$--diagonal with
amalgamation over $\D$ (see Remark \ref{newremark}).
It remains to show that $z=(x_1 x_2)(x_1 x_2)\gwia$ is free from $\D$.

We can extend the algebra $\A'$ in such a way that $\D\subset \A'$ and
sets $\D$, $\{x_1', (x_1')\gwia\}$, $\{x_2',(x_2')\gwia\}$ are free.

First of all notice that it follows from Theorem \ref{theo:schlyakhtenko}
that $\Alg(\D,x_1,x_1\gwia)$ and $x_2 x_2\gwia$ are free.
Therefore for $d_1,\dots,d_n\in\D$ we have that
$\phi[d_1 z d_2 z \cdots d_n z)\gwia]$
is uniquely determined by moments $\phi[(x_2 x_2\gwia)^k]$ and
\begin{equation}
\phi[f_1 x_1 x_1\gwia f_2 x_1 x_1\gwia \cdots f_m x_1 x_1\gwia]
\label{eq:6863}
\end{equation}
for some $f_1,\dots,f_m\in\D$.

But again Theorem \ref{theo:schlyakhtenko} gives us that $\D$ and
$x_1 x_1\gwia$ are free; therefore (\ref{eq:6863}) is uniquely
determined by moments $\phi[(x_1 x_1\gwia)^l]$ and moments
$\phi(d_{i_1} \cdots d_{i_p})$.

If we repeat the above discussion for the product
$\phi'(d_1 z' d_2 z' \cdots d_n z')$, where $z'=x_1' x_2' (x_1' x_2')\gwia$,
we see that
$$\phi(d_1 z d_2 z \cdots d_n z)=\phi'(d_1 z' d_2 z' \cdots d_n z').$$

The above result implies that $\D$ and $z$ are free and that
$z\in(\A,\phi)$ and $z'\in(\A',\phi')$ have the same distribution.
It follows that $x_1 x_2\in(\A,\phi)$ and $x_1' x_2'\in(\A',\phi')$
have the same distribution.
\end{proof}


\begin{remark}
Note that in general $x_1,x_2\in(\A,\phi)$ are not free.
\end{remark}

The above theorem gives us the upper triangular form of an operator
$x_1 x_2$ if the upper triangular forms of $x_1$ and $x_2$ are known.
Indeed, $x_1 x_2$ is an upper triangular matrix, whose diagonal elements 
are $e_k x_1 x_2 e_k=(e_k x_1 e_k) (e_k x_2 e_k)$. 

Since $x_1$ and $x_2$ are free with amalgamation over $\D$ and each of them
is $R$--diagonal with amalgamation over $\D$, we have that
in the scalar
probability space $(e_k \A e_k,\frac{1}{t_k} \phi)$ the elements 
$e_k x_1 e_k$ and $e_k x_2 e_k$ are free and each of them is $R$--diagonal.
Lemma \ref{lem:stransformataipromienspektralny} shows therefore that 
the spectral radius of $e_k x_1 x_2 e_k\in e_k \A e_k$ is equal to
$$\|e_k x_1 x_2 e_k\|_{L^2}=
\|e_k x_1 e_k\|_{L^2}  \|e_k x_2 e_k\|_{L^2}$$
which is equal to the
product of the spectral radius of $e_k x_1 e_k\in e_k \A e_k$ and
the spectral radius of $e_k x_2 e_k\in e_k \A e_k$.

We know that the spectral
radius of $e_k x_1e_k $, $e_k x_2 e_k$, and $e_k x_1 x_2 e_k$ is given by
$S_{x_1x_1\gwia }(-t)$, $S_{x_2x_2\gwia }(-t)$, and
$S_{(x_1x_2)(x_1x_2)\gwia}(-t)$, respectively, where $t=\phi(e_1+\cdots+e_k)$,
the above discussion shows that for every $-1\leq s\leq 0$ we have
$$S_{x_1 x_1\gwia} (s) S_{x_2 x_2\gwia}(s)= S_{x_1 x_2 x_2\gwia
x_1\gwia} (s). $$

\begin{corollary}
If $\mu$ and $\nu$ are probability measures on $[0,\infty)$ then
\begin{equation}S_{\mu} (s) S_{\nu} (s)=S_{\mu\boxtimes \nu} (s) \qquad
\mbox{for every } -1\leq s\leq 0, \label{eq:1989} \end{equation}
where $\boxtimes$ denotes the free multiplicative convolution \cite{VDN}.
\index{free multiplicative convolution}
\end{corollary}

\begin{remark}
Since $S$--transform is a holomorphic function it follows, that the above
equations holds for all $s$ for which it makes sense. It is also possible
to remove the assumption that $\mu$ and $\nu$ are measures on $[0,\infty)$
by the following argument: every moment of $\mu\boxtimes\nu$ is
a polynomial in moments of $\mu$ and moments of $\nu$. If we
expand $S_{\mu}(s)$, $S_{\nu}(s)$ and $S_{\mu\boxtimes\nu}(s)$ as power series
around $0$, we have that coefficients are polynomials of moments of $\mu$ and
moments of $\nu$. Since identity (\ref{eq:1989}) holds for a sufficiently many
measures it follows that appropriate polynomials on both side of the equation
are equal; therefore (\ref{eq:1989}) holds for all probabilistic measures.
\end{remark}

Thus we succeeded in finding a nice interpretation for the $S$--transform
of a nonnegative operator $a$ as appropriate spectral radii of the
$R$--diagonal operator $x$ such that $x x\gwia=a$. Furthermore this
interpretation allowed us to find a new proof of the multiplicative property
of the $S$--transform.

\section {Random matrices}
\label{chap:losowemacierze}

\subsection {Discontinuity of Brown measure}
One of the greatest difficulties connected with the Kadison--Fuglede determinant 
and Brown spectral distribution measure is that these two objects do not behave in
a continuous way.

For example, the distribution of an $n\times n$ nilpotent matrix
$$\left[ \begin{array}{ccccc} 
0 	& 0 	&\cdots	& 0 	& 0		\\
1 	& 0 	&\cdots & 0	& 0		\\
0 	& 1 	&\ddots &\vdots & \vdots	\\	
\vdots 	&\ddots &\ddots	& 0	& 0		\\
0	&\cdots	& 0	& 1	& 0
\end{array}\right]\in (\M_n,\tr_n)$$
converges, as $n$ tends to infinity, to a Haar unitary.
All of the above matrices have the determinant equal to $0$, while the Haar unitary
has the Fuglede--Kadison determinant equal to $1$; all of the above matrices
have the Brown measure equal to $\delta_0$, while the Brown measure of the Haar
unitary is the uniform measure on the unit circle $\{z\in\C:|z|=1\}$.

We would like to discuss the following problem: suppose that a Kolmogorov
probability space $(\Omega,{\mathcal{B}},P)$ and a
sequence of
random matrices $m_n\in {\mathcal{L}}^{\infty -}\otimes \M_n$ are given.
We can regard the set of $n\times n$ random matrices with all moments finite
as a noncommutative probability space with a state $\phi_n (z)=\E \tr_n z$.
Suppose that the distribution of $m_n$ converges to a distribution of a certain
element $x\in (\A,\phi)$, where $\A$ is a von Neumann algebra.
Our goal is to determine the empirical distribution
of eigenvalues of $m_n$ or---in other words---the Brown measure of $m_n$.
In the light of the previous discussion we should not expect that the sequence
of Brown measures of $m_n$ will always converge to the Brown measure of
$x\in(\A,\phi)$.

Surprisingly, in many known cases when we consider a ``reasonable'' sequence
of random matrices the sequence of Brown measures converges to the Brown measure
of the limit (cf \cite{BL}). In these cases, however, the convergence was proved
in this way, that the distribution of eigenvalues of $m_n$
was calculated by ad hoc methods and nearly by
an accident it turned out to converge to the Brown measure of $x$.
Therefore one of the
most interesting problems in the theory of random matrices is to relate the
asymptotic distribution of eigenvalues of a sequence of random matrices with
the Brown measure of the limit---for example to clarify what does it mean that
a sequence of random matrices is ``reasonable''.

In this Section we present an interesting direction of the future research,
which possibly may give light into this subject.

\subsection {Asymptotical freeness}
It was observed by Voiculescu that many random matrix models provide very
natural examples of noncommutative random variables which are, informally
speaking, free when the size of the matrices tends to infinity. He
introduced the following definition.
\begin{definition}
Let $(\Omega,{\mathcal{B}},P)$ be a probability space and
let $(a_{1,n})$ and $(a_{2,n})$ be sequences of random matrices,
$a_{1,n},a_{2,n}\in{\mathcal{L}}^{\infty -}(\Omega,\M_n(\C))$.

Let $\A$ be the free $\star$--algebra generated by elements $a_1$ and $a_2$
and let $\phi:\A\rightarrow\C$ be a state defined by
$$\phi(a_{i_1}^{s_1} \cdots a_{i_k}^{s_k})=
\lim_{n\rightarrow\infty}
{\mathbb{E}}\big[ {\rm tr}_n a_{i_1,n}^{s_1} \cdots a_{i_k}^{s_k}\big]$$
for every $k\in\N$, $i_1,\dots,i_k\in\{1,2\}$ and
$s_1,\dots,s_k\in\{1,\star\}$, where ${\rm tr}_n=\frac{1}{n} {\rm Tr}$
denotes the normalized trace on $\M_n(\C)$.

We say that the sequences $(a_{1,n})$ and $(a_{2,n})$ are asymptotically free
\index{asymptotical freeness}
if the above limit always makes sense and if $\{a_1,a_1\gwia\}$ and
$\{a_2,a_2\gwia\}$ are free.
\end{definition}

\subsection {Asymptotical freeness conjecture}

Let $(D_n)$ be a sequence of constant hermitian matrices,
$D_n\in\M_n(\C)$, such that the distribution of the sequence $D_n$
converges (i.e.\ for each $k\in\N$ the limit
$\lim_{n\rightarrow\infty} {\mathbb{E}}\ {\rm tr}_n D_n^k$ exists) and
such that the operator norms are uniformly bounded: $\|D_n\|<C$.

For every $n\in\N$ let $V_n$ be a classical random variable with values in
the group of unitary matrices $U(n)$, which is Haar distributed
on this group.

We define a random matrix $x_n=V_n D_n$---this is the
so--called $U(n)$--biinvariant ensemble.
It follows from results of Voiculescu that the $\star$--distribution of the
random variable $x_n\in(\M_n(\C),{\rm tr}_n)$ converges in moments
to a distribution of an $R$--diagonal element as $n$ tends to infinity and
that for each $k$ the sequences $(\sqrt[2k]{ (x_n\gwia)^k x_n^k})_{n\in\N}$
and $x_n x_n\gwia$ are asymptotically free. We hope therefore that the
following conjecture is true.
\begin{conjecture}
The sequences $(\lim_{k\rightarrow\infty} \sqrt[2k]{(x_n\gwia)^k
x_n^k})_{n\in\N}$ and $(x_n x_n\gwia)_{n\in\N}$ are asymptotically free.
\end{conjecture}

For every $n$ and $\omega\in\Omega$ we can write the matrix $x_n(\omega)$
in the upper triangular form $y_n(\omega)$, cf (\ref{eq:utm});
the question is to determine the joint distribution of the random variables
$\lambda_1,\dots,\lambda_n$ and $a_{ij}$ ($1\leq i<j\leq n$).

The above conjecture would give us a partial solution to the above
problem; for fixed $t_1,\dots,t_N>0$ such that $t_1+\cdots+t_N=1$
we consider sequences $k_0(n),\dots,k_N(n)$ such that for $n$
large enough we have $0=k_0(n)<k_1(n)<\cdots<k_N(n)=N$ are integer
numbers and for each $1\leq m\leq N$ we have
$$\lim_{n\rightarrow\infty} \frac{k_m(n)}{n}=t_1+\cdots+t_{m}.$$
We consider the corresponding orthogonal projections
$e_1(n),\dots,e_N(n)\in\M_n(\C)$ given by
$$e_i(n)={\rm diag}(\underbrace{0,\dots,0}_{k_{i-1}(n)},
\underbrace{1,\dots,1}_{k_i(n)-k_{i-1}(n)},
\underbrace{0,\dots,0}_{n-k_i(n)}).$$

Using these projections $e_i$ we can divide the triangular matrix
$y_n(\omega)$ into submatrices:
$$y_n(\omega)=\left[ \begin{array}{ccccc} e_1 y_n(\omega) e_1 &
e_1 y_n(\omega) e_2 & \cdots & e_1 y_n(\omega) e_{n-1}& e_1 y_n(\omega) e_n
\\ 0 & e_2 y_n(\omega) e_2 &\cdots & e_2 y_n(\omega) e_{n-1}& e_2
y_n(\omega) e_n \\ \vdots& & \ddots & \vdots & \vdots \\
& & & e_{n-1} y_n(\omega) e_{n-1} & e_{n-1} y_n(\omega) e_n \\
0& & \cdots & 0 & e_n y_n(\omega) e_n \end{array} \right]. $$

{} From the conjecture it would follow that $\{e_1(n),\dots,e_N(n)\}$ and
$y_n(\omega) y_n\gwia(\omega)$ are asymptotically free; therefore the joint
distribution of operators $e_i y_n(\omega) e_j$ converges to the
distribution of operators $x_{ij}$ constructed in Theorem
\ref{theo:utmform}.

It would follow that the Brown measures of the considered random matrices
converge weakly to the Brown measure of the limit distribution.

\end{document}